\documentclass[11pt,final]{amsart}
\usepackage[utf8]{inputenc}
\usepackage[T1]{fontenc}
\usepackage[letterpaper,centering]{geometry}
\usepackage{lmodern}
\usepackage{eucal}
\usepackage{mathrsfs}
\usepackage{amsmath,amssymb,amsfonts,amscd}
\usepackage{microtype,xcolor,enumitem,booktabs,needspace}
\definecolor{darkblue}{rgb}{0,0,0.3}
\definecolor{darkgreen}{rgb}{0,0.4,0}
\usepackage[colorlinks=true,linkcolor=darkblue,citecolor=darkgreen,
urlcolor=darkblue,unicode,pdfborder={0 0 0}]{hyperref}
\hypersetup{pdftitle={Brauer groups of generalised Markoff surfaces},
	pdfauthor={Zheshuai Tian and Dasheng Wei},
	pdfsubject={Brauer groups, Markoff surfaces, conic bundles and Galois cohomology}}
\setlist[enumerate]{label={\upshape(\arabic*)},topsep=.7ex,leftmargin=*}
\setlist[itemize]{leftmargin=*}
\theoremstyle{plain}
\newtheorem{thm}{Theorem}[subsection]
\renewcommand{\thethm}{%
	\ifnum\value{subsection}=0
	\thesection.\arabic{thm}%
	\else
	\thesubsection.\arabic{thm}%
	\fi
}
\newtheorem{lem}[thm]{Lemma}
\newtheorem{prop}[thm]{Proposition}
\newtheorem{cor}[thm]{Corollary}
\theoremstyle{definition}
\newtheorem{defn}[thm]{Definition}
\newtheorem{rmk}[thm]{Remark}
\newtheorem{example}[thm]{Example}
\numberwithin{equation}{section}
\newcommand{\PP}{\mathbf P}
\newcommand{\Aff}{\mathbf A}

\newcommand{\FF}{\mathbf F}
\newcommand{\ZZ}{\mathbf Z}
\newcommand{\QQ}{\mathbf Q}
\newcommand{\kbar}{\overline{k}}
\newcommand{\Gm}{\mathbf G_{\mathrm m}}

\DeclareMathOperator{\Br}{Br}
\DeclareMathOperator{\Pic}{Pic}
\DeclareMathOperator{\Gal}{Gal}
\DeclareMathOperator{\Hom}{Hom}

\title[Brauer groups of Markoff surfaces]
{Brauer groups of generalised Markoff surfaces}

\author{Zheshuai Tian}
\address{Academy of Mathematics and Systems Science, CAS, Beijing 100190,
	P.\ R.\ China}
\email{tianzheshuai@amss.ac.cn}
\author{Dasheng Wei}
\address{Academy of Mathematics and Systems Science, CAS, Beijing 100190,
	P.\ R.\ China {\normalfont\itshape and} School of Mathematical Sciences,
	University of CAS, Beijing 100049, P.\ R.\ China}
\email{dshwei@amss.ac.cn}
\date{}

\begin{document}
	
	\begin{abstract}
		Let $k$ be a field of characteristic zero with $a,b,c\in k^\times$ 
		and let $m\in k\setminus\{0,4abc\}$.  We compute
		the Brauer group of the generalised Markoff surface
		$ax^2+by^2+cz^2-xyz=m$.  We also construct several counterexamples to the
		integral Hasse principle, which are explained by the Brauer--Manin obstruction.
	\end{abstract}
	
	\maketitle
	
	\section{Introduction}
	
	Markoff introduced the equation
	\[
	x^2+y^2+z^2=3xyz
	\]
	\cite{Mar79}, whose positive integral solutions form an orbit under
	the action of a discrete group \cite{BU04,BGS16}.
	With the change of variables $(X,Y,Z)=(3x,3y,3z)$,
	the equation becomes $X^2+Y^2+Z^2-XYZ=0$, 
	which defines a singular Markoff surface.
	
	Let $k$ be a field of characteristic zero with $a,b,c\in k^\times$ 
	and let $m\in k\setminus\{0,4abc\}$.   
	We study the smooth projective cubic surface
	\begin{equation}\label{eq:projective}
		X_m:\quad (ax^2+by^2+cz^2)w-xyz=mw^3\subset\PP^3_k
	\end{equation}
	and its affine part, the generalised Markoff surface
	\begin{equation}\label{eq:affine}
		U_m:\quad ax^2+by^2+cz^2-xyz=m.
	\end{equation}
	The boundary is the split triangle
	\begin{equation}\label{eq:boundary}
		\mathscr D=D_x+D_y+D_z,\qquad
		D_x=\{x=w=0\},\quad D_y=\{y=w=0\},\quad D_z=\{z=w=0\}.
	\end{equation}

	For ordinary Markoff surfaces $x^2+y^2+z^2-xyz=m$,
	Ghosh and Sarnak \cite{GS22} proved that the integral Hasse principle
	holds for almost all integer parameters and exhibited infinitely many
	failures. Colliot-Th\'el\`ene, Xu and the second author
	\cite{CTWX20} and Loughran and Mitankin \cite{LM21} computed
	Brauer groups and studied the integral Brauer--Manin obstruction.
	S.~Chen \cite{SChen21} computed the algebraic Brauer group
	of coefficient twists $ax^2+y^2+z^2-xyz=m$ and constructed
	examples with order-four algebraic Brauer classes.
	Dao \cite{Dao24} studied Brauer groups and integral
	local--global principles for other Markoff-type cubic surfaces.
	For the general family \eqref{eq:affine}, 
	we compute both the Brauer group $\Br(U_m)$ and the unramified Brauer group $\Br(X_m)$.
	
	Put
	\begin{equation}\label{eq:Delta}
		\Delta=m-4abc.
	\end{equation} 
	We begin with our computation on the unramified Brauer group, 
	which is proved in Theorem~\ref{proj:classification}.
	\begin{thm} \label{thm:main-projective}
	Let $k$ be a field of characteristic zero with $a,b,c\in k^\times$ 
	and let $m\in k\setminus\{0,4abc\}$.
		One has $\Br(X_m)=\Br_1(X_m)$ and
		\begin{equation}\label{eq:intro-projective}
			\Br(X_m)/\Br(k)\simeq
			\begin{cases}
				0,&\{a\Delta,b\Delta,c\Delta,m\Delta\}\cap k^{\times2}
				\ne\varnothing,\\
				(\ZZ/2)^2,&ab,ac,am\in k^{\times2}
				\text{ and }a\Delta\notin k^{\times2},\\
				\ZZ/2,&\text{otherwise}.
			\end{cases}
		\end{equation}
	\end{thm}

	We next compute the algebraic part $\Br_1(U_m)/\Br(k)$. Write
	\begin{equation}\label{eq:intro-characters}
		\chi_1=[a\Delta],\qquad\chi_2=[b\Delta],\qquad
		\chi_3=[c\Delta],\qquad\chi_4=[m\Delta]
		\quad\text{in }k^\times/k^{\times2}.
	\end{equation}
	These four square classes determine the Galois action on \(\Pic(\overline U_m)\), which gives the following
	classification. 
	
	\begin{thm} \label{thm:main-algebraic}
		Let $k$ be a field of characteristic zero with $a,b,c\in k^\times$ 
		and let $m\in k\setminus\{0,4abc\}$.
		Let $N$ be the number of nonzero characters in
		\eqref{eq:intro-characters} and put
		$\rho=\dim_{\FF_2}\langle\chi_1,\ldots,\chi_4\rangle$.
		Then
		\begin{equation}\label{eq:intro-algebraic}
			\Br_1(U_m)/\Br(k)\simeq
			\begin{cases}
				(\ZZ/2)^{N-\rho},&N<4,\\
				(\ZZ/2)^{5-\rho},&N=4\text{ and every relation has even weight},\\
				\ZZ/4\oplus(\ZZ/2)^{3-\rho},
				&N=4\text{ and an odd-weight relation exists}.
			\end{cases}
		\end{equation}
		Here a relation is a vector
		$\epsilon=(\epsilon_1,\ldots,\epsilon_4)\in\FF_2^4$ such that
		\[
		\epsilon_1\chi_1+\epsilon_2\chi_2+
		\epsilon_3\chi_3+\epsilon_4\chi_4=0.
		\]
		Its weight is defined by $\#\{i:\epsilon_i=1\}$ and an odd-weight relation
		is one whose weight is odd.
	\end{thm}
	
	For the transcendental part $\Br(U_m)/\Br_1(U_m)$ over a number field, put
	\[
	M=\frac{m}{abc},\qquad D=\frac{\Delta}{abc}=M-4,
	\]
	and choose square roots in a fixed algebraic closure of $k$.
	Define
	\[
	\begin{gathered}
		K_1=k(\sqrt{a\Delta}),\qquad K_2=k(\sqrt{b\Delta}),\qquad
		K_3=k(\sqrt{c\Delta}),\\
		L=k(\sqrt{a\Delta},\sqrt{b\Delta},\sqrt{c\Delta},
		\sqrt{m\Delta}),\qquad
		\eta=\frac{\sqrt M-\sqrt D}{2}.
	\end{gathered}
	\]
	Both square roots defining $\eta$ belong to $L$, and
	$\eta\in L^\times$ because $M-D=4$.
	For a number field $F$, let $w_F=\#\mu(F)$.
	
	\begin{thm}
		\label{thm:main-transcendental}
		Let $k$ be a number field with $a,b,c\in k^\times$ and let $m\in k\setminus\{0,4abc\}$.  With the notation above, put
		\[
		E=\gcd\left(w_k,\frac{w_{K_1}}2,
		\frac{w_{K_2}}2,\frac{w_{K_3}}2\right),
		\qquad h=[L(\sqrt[E]{\eta}):L].
		\]
		Then $h\mid E$ and
		\begin{equation}\label{eq:intro-transcendental}
			\frac{\Br(U_m)}{\Br_1(U_m)}\simeq\ZZ/n_0,
			\qquad n_0=\frac{E}{h}.
		\end{equation}
		%For every positive integer $n$, we have $n\mid n_0$ if and only if
		%$n\mid E$ and $\eta\in L^{\times n}$.  In particular, $n_0\mid w_k$.
		The result is independent of the choices of roots.
	\end{thm}
	
	Theorem~\ref{thm:complete-affine-brauer} determines the
	extension of the transcendental part by the algebraic part.
	Together with the splitting of the constant subgroup $\Br(k)$,
	this gives the full Brauer group $\Br(U_m)$.

	In particular, if $k$ is a number field and has a real embedding (for example, $k=\QQ$), then
	$\mu(k)=\{\pm1\}$. Hence
	the transcendental part is either $0$ or $\ZZ/2$.
	In this case, square conditions in $k$ determine
	the full Brauer group and whether the extension splits.
	
	\begin{thm}
		\label{thm:main-real}
		Let $k$ be a number field with a real embedding
		$v_0:k\hookrightarrow\mathbf R$, let $a,b,c\in k^\times$, and
		suppose $m\ne0,4abc$.  Set
		\[
		A=\Br_1(U_m)/\Br(k),\qquad B=\Br(U_m)/\Br(k),
		\]
		where $A$ is given by Theorem~\ref{thm:main-algebraic}.
		If at least one of $-a\Delta,-b\Delta,-c\Delta$ is not a square
		in $k$, then $B=A$.
		Otherwise $-D\in k^{\times2}$; choose $d\in k$ with
		$d^2=-D$ and $v_0(d)>0$, and put
		\[
		p=1+\frac d2,\qquad q=1-\frac d2.
		\]
		Then $A$ is killed by $2$, and
		\begin{equation}\label{eq:intro-real-full}
			B\simeq
			\begin{cases}
				A\oplus\ZZ/2,
				&p\in k^{\times2}\text{ or }q\in k^{\times2},\\
				\ZZ/4\oplus(\ZZ/2)^2,
				&-q\in k^{\times2}\text{ and }p\notin k^{\times2},\\
				A,&\text{otherwise}.
			\end{cases}
		\end{equation}
		In the first two cases, $B/A\simeq\ZZ/2$; its extension by $A$
		is split in the first case and nonsplit in the second.
		%In all cases $\Br(U_m)\simeq\Br(k)\oplus B$ noncanonically.
	\end{thm}
	
	We also apply these computations to integral points.
	In Section~\ref{sec:integral-failures}, we construct two explicit
	families over $\QQ$ with real points and $\ZZ_p$-points for
	every prime $p$, but no integral points.
	For each case, a quaternion class gives a Brauer--Manin
	obstruction. The constructions follow the reciprocity
	approach of \cite{GS22,CTWX20,SChen21}.

	The paper is organized as follows. Section~\ref{sec:prelim} develops
	the geometric and descent preliminaries, and Section~\ref{proj:section}
	computes the unramified Brauer group. Section~\ref{sec:algebraic}
	computes the Picard lattice and the algebraic part, while
	Section~\ref{sec:complete-affine} determines the full affine Brauer
	group. Section~\ref{sec:consequences} specializes the results to
	number fields with a real embedding and to ordinary
	Markoff surfaces. Finally,
	Section~\ref{sec:integral-failures} treats the integral Hasse failures.

	\subsection*{Notation and conventions}
	
	Fix an algebraic closure $\kbar$ of $k$ and put
	$G_k=\Gal(\kbar/k)$. 
	Cohomology of varieties means \'etale cohomology, and cohomology of
	fields and Galois modules means continuous Galois cohomology.
	We use the cohomological Brauer group
	$\Br(V)=H^2(V,\Gm)_{\mathrm{tors}}$ and write
	\[
	\Br_1(V)=\ker[\Br(V)\longrightarrow\Br(\overline V)].
	\]
	The notation $(f,g)$ denotes a quaternion class.  Square classes are
	written additively, with $[q]=0$ meaning $q\in k^{\times2}$.
	For a torsion group $A$, $A[n]$ denotes its subgroup killed by $n$.
	The Tate twist $\ZZ/n(-1)$ means $\Hom(\mu_n,\ZZ/n)$.
	The constant Brauer groups inject in all cases considered here, so
	we identify $\Br(k)$ with its image. We refer to
	$\Br_1(V)/\Br(k)$ and $\Br(V)/\Br_1(V)$ as the algebraic and
	transcendental parts, respectively.
	
	\section{Geometry and descent}\label{sec:prelim}
	
	\begin{lem}\label{lem:smooth-points}
		The surface $X_m$ is smooth and geometrically rational, and
		$U_m(k)\ne\varnothing$. It follows that both
		$\Br(k)\to\Br(X_m)$ and $\Br(k)\to\Br(U_m)$ are split injections,
		and $\Br(\overline X_m)=0$.
	\end{lem}
	
	\begin{proof}
		On the boundary, the derivatives with respect to $x,y,z$ are
		$-yz,-xz,-xy$.  They vanish simultaneously only at the three
		vertices of the triangle, where the derivative with respect to $w$
		takes the values $a$, $b$, and $c$, respectively. Thus there are  no singular points at
		infinity.  In the affine chart, a singular point satisfies
		\[
		2ax=yz,\qquad 2by=xz,\qquad 2cz=xy.
		\]
		If a coordinate vanishes, all three vanish, which forces $m=0$.
		Otherwise these equations give
		$x^2=4bc$, $y^2=4ac$, $z^2=4ab$, and $xyz=8abc$, so the
		surface equation forces $m=4abc$.  Both possibilities are excluded.
		
		Every smooth cubic surface over an algebraically closed field is
		rational; see, for example, \cite{Poonen}.  The boundary vertex
		$[1:0:0:0]$ is a $k$-point.  A smooth cubic surface with a
		$k$-point is $k$-unirational \cite[Remark~9.4.11]{Poonen}.
		Since $k$ is infinite, unirationality implies that $X_m(k)$ is Zariski dense in $X_m$;
		hence the nonempty open subset $U_m$ has a $k$-point.
		Evaluation at these points splits the inclusion of the constant Brauer group.
		Finally, geometric rationality gives $\Br(\overline X_m)=0$.
	\end{proof}
	
	\begin{lem}\label{lem:descent}
		Every invertible regular function on $\overline U_m$ is constant,
		and $\Pic(\overline U_m)$ is a free abelian group of rank four.  For
		$V=X_m$ or $U_m$,  one has
		\[
		\Pic(V)=\Pic(\overline V)^{G_k},\qquad
		\Br_1(V)/\Br(k)=H^1(k,\Pic(\overline V)).
		\]
		For every $n\ge1$,
		\begin{equation}\label{eq:prelim-HS}
			H^2(U_m,\mu_n)/H^2(k,\mu_n)
			\simeq H^2(\overline U_m,\mu_n)^{G_k}.
		\end{equation}
	\end{lem}
	
	\begin{proof}
		The assertions about invertible regular functions and the Picard
		group of $\overline U_m$ follow from the calculation in
		Section~\ref{subsec:geometric-configuration}.
		
		For the assertions about the Picard and algebraic Brauer groups,
		apply the Hochschild--Serre spectral sequence to $\Gm$.  A rational point makes the pullback
		$H^j(k,\Gm)\to H^j(V,\Gm)$ injective in every degree, so the
		usual obstruction maps vanish.  This gives the stated identifications.
		For \eqref{eq:prelim-HS}, the Kummer sequence gives
		$H^1(\overline U_m,\mu_n)=0$ because the invertible regular
		functions are constant and the geometric Picard group is torsion-free.  The only possible
		differential out of $H^2(\overline U_m,\mu_n)^{G_k}$ is then the
		$d_3$ into $H^3(k,\mu_n)$; it vanishes because the rational point
		also splits this pullback.  The edge sequence is
		\eqref{eq:prelim-HS}.  We use these standard descent and Kummer
		constructions in the form explained in \cite{Milne}.
	\end{proof}

	\section{The Brauer group of the projective surface}\label{proj:section}
	
	Throughout this section put $\Delta=m-4abc$ and let $t=x/w$.
    We compute classes unramified along all prime divisors, including the three
    lines at infinity.
	
	\begin{thm}\label{proj:classification}
		The Brauer group of $X_m$, modulo constant classes, is
		\[
		\Br(X_m)/\Br(k)\cong
		\begin{cases}
			0,&\text{if at least one of }a\Delta,b\Delta,c\Delta,m\Delta
			\text{ belongs to }k^{*2},\\[2pt]
			(\ZZ/2)^2,&\text{if }ab,ac,am\in k^{*2}
			\text{ and }a\Delta\notin k^{*2},\\[2pt]
			\ZZ/2,&\text{otherwise}.
		\end{cases}
		\]
		In every nontrivial case the quaternion class
		\[
		B=(t^2-4bc,b\Delta)
		\]
		belongs to $\Br(X_m)$ and is nonconstant.  It generates this quotient
		whenever the quotient is cyclic.  In the second case, choose
		$u,v\in k^*$ with $u^2=bc$ and $v^2=m/a$; these exist since
		$bc=(ab)(ac)/a^2$ and $m/a=(am)/a^2$ are squares.  Then the classes
		\[
		B=(t^2-4bc,b\Delta),\qquad
		C=((t-2u)(t-v),b\Delta)
		\]
		generate $\Br(X_m)/\Br(k)$.
	\end{thm}
	
	We first describe the conic pencil and its splitting characters.
	Let $L_x$ be the line $x=w=0$ on $X_m$.  The planes $x=tw$
	through $L_x$, after removal of their fixed component $L_x$, define
	the conic-bundle morphism
	\[
	f:X_m\longrightarrow\PP^1_k.
	\]
	On $w\ne0$ this morphism is given by the function $t=x/w$.
	Its finite fibres are the residual conics
	\begin{equation}\label{proj:conics}
		by^2+cz^2-tyz+(at^2-m)w^2=0.
	\end{equation}
	The fibre at infinity is the union of the two $k$-lines
	$y=w=0$ and $z=w=0$.  In particular, this fibre is split.
	The two quadratics
	\[
	F_1(t)=t^2-4bc,\qquad F_2(t)=t^2-m/a
	\]
	are separable and have disjoint zero sets.  Their four geometric
	zeros correspond to all the finite degenerate fibres of the pencil.
	
	\begin{lem}\label{proj:splitting}
		Let $Q$ be a closed point of $\PP^1 \backslash \{\infty\}=\Aff^1_k$.  If $F_1(Q)=0$,
		the two geometric components of $f^{-1}(Q)$ are defined over
		$k(Q)(\sqrt{b\Delta})$.  If $F_2(Q)=0$, they are defined over
		$k(Q)(\sqrt{a\Delta})$.  The corresponding generic conic has
		Brauer class
		\begin{equation}\label{proj:conic-class}
			A=(F_1,b(m-at^2))=(F_1,b\Delta)+(F_2,a\Delta)
			\quad\text{in }\Br(k(t)).
		\end{equation}
	\end{lem}
	
	\begin{proof}
		Multiplying \eqref{proj:conics} by $4b$ and completing the square gives
		\[
		(2by-tz)^2=F_1z^2+4b(m-at^2)w^2.
		\]
		At a zero of $F_1$, this becomes
		$(2by-tz)^2=4b\Delta w^2$, proving the first assertion.
		At a zero of $F_2$, the equation is binary in $y,z$, and its
		discriminant is $F_1(Q)=\Delta/a$, whose square class is
		$a\Delta$.  This establishes the second assertion and yields the first expression
		for $A$.
		
		For the final identity, observe that $F_1-F_2=\Delta/a$.
		The quaternion identity
		\[
		(r,r-s)=(r,-s)+(s,r-s)
		\]
		follows from the Steinberg relation $(u,1-u)=0$ for $u\ne0,1$
		(see \cite[Example~1.3.4, p.~8]{GS17}), with $u$ being replaced by $s/r$.
		Taking $r=F_1$ and $s=\Delta/a$, and adding $(F_1,-ab)$ to
		both sides gives
		\[
		\begin{aligned}
			A=(F_1,-abF_2)
			&=(F_1,(-ab)(-\Delta/a))+(a\Delta,F_2)\\
			&=(F_1,b\Delta)+(F_2,a\Delta),
		\end{aligned}
		\]
		where $m-at^2=-aF_2$.  This uses the bilinearity and symmetry of
		quaternion symbols and proves \eqref{proj:conic-class}.
	\end{proof}
	
	Let $\mathcal N$ be the set of closed points whose fibres are nonsplit degenerate. 
	 For $Q\in\mathcal N$ write
	\[
	d_Q=\begin{cases}b\Delta,&F_1(Q)=0,\\
		a\Delta,&F_2(Q)=0.
	\end{cases}
	\]
	By definition, $d_Q$ represents a nonzero quadratic character of $k(Q)$.
	Identify quadratic characters with square classes and set
	\begin{equation}\label{proj:residue-kernel}
		R=\ker\left[\FF_2^{\mathcal N}\longrightarrow k^*/k^{*2},\quad
		(\epsilon_Q)_Q\longmapsto
		\prod_{Q\in\mathcal N}N_{k(Q)/k}(d_Q)^{\epsilon_Q}\right].
	\end{equation}
	
	\begin{prop}\label{proj:residue-description}
		If $\mathcal N\ne\varnothing$, there is an isomorphism
		\[
		\Br(X_m)/\Br(k)\cong R/\langle\boldsymbol 1\rangle,
		\qquad \boldsymbol 1=(1)_{Q\in\mathcal N}.
		\]
		If $\mathcal N=\varnothing$, then $\Br(X_m)=\Br(k)$.
	\end{prop}
	
	\begin{proof}
		By \cite[Lemma~11.3.2]{CTSk21}, contracting the split fibre
		components gives a smooth projective conic bundle
		$X'\to\PP^1_k$ with reduced irreducible fibres and
		$\Br(X')=\Br(X_m)$.  Its singular fibres are indexed by
		$\mathcal N$, with splitting characters $d_Q$; its fibre at
		infinity is smooth and split.
		If $\mathcal N\ne\varnothing$, the generic conic class $A$ is
		ramified, and \cite[Proposition~11.3.4]{CTSk21} identifies
		$\Br(X_m)/\Br(k)$ with $R/\langle\boldsymbol1\rangle$, since
		the residue vector of $A$ is $\boldsymbol1$.  If $\mathcal N=\varnothing$, then $A$ is
		constant and specializes to zero at infinity, hence $A=0$;
		thus $X'$ is $k$-rational and $\Br(X_m)=\Br(k)$.
		
		Since $d_Q\in k^*$, \cite[Corollary~11.3.5]{CTSk21} gives
		representatives
		\[
		\beta_\epsilon=\sum_{Q\in\mathcal N}
		\epsilon_Q(p_Q(t),d_Q),\qquad \epsilon\in R,
		\]
		where $p_Q$ is the monic irreducible polynomial of $Q$.
		Their pullbacks represent the classes in the quotient.
	\end{proof}
	
	\begin{proof}[Proof of Theorem~\ref{proj:classification}]
		The fibres over the zeros of $F_1$ are split exactly when
		$b\Delta$ or $c\Delta$ is a square in $k$.  Indeed, if $bc$
		is a square, both zeros are rational and the criterion is simply
		$b\Delta\in k^{*2}$.  Otherwise their residue field is
		$k(\sqrt{bc})$, and a scalar from $k^*$ becomes a square in
		this quadratic field exactly when its square class is $1$ or $bc$.
		Similarly the fibres over the zeros of $F_2$ are split exactly
		when $a\Delta$ or $m\Delta$ is a square.
		
		If the fibres in both pairs are split, then $\mathcal N=\varnothing$, so
		Proposition~\ref{proj:residue-description} gives a zero quotient.
		If the fibres in exactly one pair are split, then all nonsplit fibres belong to the
		other pair.  That pair contributes either one quadratic closed
		point or two rational points with the same character; in either
		case $R=\langle\boldsymbol1\rangle$, so the quotient is again zero.
		The first case of the theorem follows.
		
		Suppose now that neither pair is split.  An irreducible quadratic
		contributes one coordinate to $R$ and no corestriction condition,
		since $N_{k(Q)/k}(d_Q)=d_Q^2$.  A rational pair contributes two
		coordinates and one condition given by their common nonzero
		character.  If both quadratics are irreducible, $R$ has dimension
		two; if exactly one splits into rational points, $R$ again has dimension
		two.  After quotienting by $\boldsymbol1$, the result in
		each case is $\ZZ/2$.
		
		If both pairs are rational, the corestriction map has rank one
		when $[b\Delta]=[a\Delta]$ and rank two otherwise.  Its kernel
		therefore has dimension three or two, respectively.  Equality of
		these two nonzero characters means $ab\in k^{*2}$.  Together
		with $bc,m/a\in k^{*2}$, this is equivalent to
		$ab,ac,am\in k^{*2}$, yielding the asserted rank-two case.
		
		By \cite[Corollary~11.3.5]{CTSk21}, $B$ corresponds to the
		full residue vector on the first pair and zero on the second.
		When neither pair is split, this vector is neither $0$ nor
		$\boldsymbol1$.  Hence $B$ has nonzero image and generates the
		quotient in the rank-one case.
		In the rank-two case, ordering the fibres as $2u,-2u,v,-v$,
		the vectors of $B,C$ are $(1,1,0,0)$ and $(1,0,1,0)$,
		which form a basis modulo $\langle(1,1,1,1)\rangle$.
		The displayed classes extend and generate.
	\end{proof}
	
	\begin{rmk}\label{proj:symmetry}
		Although the calculation uses the $x$-pencil, the result is symmetric
		in $a,b,c$. The case $(\ZZ/2)^2$ occurs exactly when $a,b,c,m$ have the
		same square class and $a\Delta$ is nonsquare. In this case, choose
		compatible square roots of $bc$, $ac$, and $ab$ and set
		\[
		X=x/\sqrt{bc},\qquad Y=y/\sqrt{ac},\qquad Z=z/\sqrt{ab}.
		\]
		Then
		\[
		X^2+Y^2+Z^2-XYZ=\frac{m}{abc},
		\]
		so this case reduces over $k$ to the usual Markoff calculation.
		For general coefficients, this change of variables need not be defined
		over $k$.
	\end{rmk}
	
	\section{The Picard lattice and the algebraic Brauer group}
	\label{sec:algebraic}
	
	Throughout this section, put $\Delta=m-4abc$ and assume
	$abcm\Delta\ne0$.  For $q\in k^\times$, write $[q]$ for its class in
	$k^\times/k^{\times2}$.  For convenience, we identify this group with
	$H^1(k,\ZZ/2)$, and write its operation additively.
	
	\subsection{The geometric configuration}
	\label{subsec:geometric-configuration}
	
	Choose $q_a,q_b,q_c,q_m\in\kbar$ with
	\[
	q_a^2=a\Delta,\qquad q_b^2=b\Delta,\qquad
	q_c^2=c\Delta,\qquad q_m^2=m\Delta,
	\]
	and put
	\[
	r_1=\frac{q_b q_c}{\Delta},\qquad
	r_2=\frac{q_a q_c}{\Delta},\qquad
	r_3=\frac{q_a q_b}{\Delta}.
	\]
	Then $r_1^2=bc$, $r_2^2=ac$, $r_3^2=ab$ and
	$r_1r_2r_3=abc$.  Set
	\begin{equation}\label{eq:alg-coordinate-change}
		X=x/r_1,\qquad Y=y/r_2,\qquad Z=z/r_3,\qquad
		M=m/(abc),\qquad D=\Delta/(abc).
	\end{equation}
	Dividing the defining cubic by $abc$ gives
	\begin{equation}\label{eq:alg-normalization}
		(X^2+Y^2+Z^2)w-XYZ=Mw^3,\qquad D=M-4.
	\end{equation}
	The diagonal projective automorphism
	$[x:y:z:w]\mapsto[x/r_1:y/r_2:z/r_3:w]$ identifies
	$\overline X_m$ with the ordinary Markoff surface of parameter
	$M=m/(abc)\notin\{0,4\}$.  It identifies their affine parts
	and carries the ordered boundary triangle to the corresponding
	triangle on the ordinary surface.  We may therefore apply the
	geometric results of \cite{CTWX20} after this change of coordinates.
	
	Write $L_1=D_x$, $L_2=D_y$ and $L_3=D_z$ for the boundary
	lines.  The geometric Picard lattice of a smooth cubic surface is
	unimodular of rank seven.
	
	The divisor sequence (see also \cite[Proposition~2.2]{CTWX20}) is
	\[
	0\longrightarrow\kbar[U_m]^\times/\kbar^\times
	\longrightarrow\bigoplus_{j=1}^3\ZZ[L_j]
	\longrightarrow\Pic(\overline X_m)
	\longrightarrow\Pic(\overline U_m)\longrightarrow0.
	\]
	The boundary intersection matrix has diagonal entries $-1$ and
	off-diagonal entries $1$, with determinant $4$.  The boundary
	classes are therefore independent, and every invertible regular
	function on $\overline U_m$ is constant.  It follows that
	\begin{equation}\label{eq:alg-pic-quotient}
		\kbar[U_m]^\times=\kbar^\times,\qquad
		\Pic(\overline U_m)=
		\Pic(\overline X_m)/\langle L_1,L_2,L_3\rangle.
	\end{equation}
	By \cite[Proposition~2.2]{CTWX20}, this quotient is torsion-free,
	hence free of rank four.  In particular, the boundary span is primitive.
	
	\begin{prop}\label{prop:alg-d4}
		There is an isomorphism of Galois lattices
		\[
		\Pic(\overline U_m)\simeq
		\mathcal L:=\{(z_1,z_2,z_3,z_4)\in\ZZ^4:
		z_1+z_2+z_3+z_4\equiv0\pmod2\},
		\]
		under which the Galois group acts diagonally through the quadratic characters
		\begin{equation}\label{eq:alg-characters}
			\chi_1=[a\Delta],\quad \chi_2=[b\Delta],\quad
			\chi_3=[c\Delta],\quad \chi_4=[m\Delta].
		\end{equation}
		Its splitting field is therefore
		\[
		K=k(\sqrt{a\Delta},\sqrt{b\Delta},\sqrt{c\Delta},\sqrt{m\Delta})
		=k(\sqrt{ab},\sqrt{ac},\sqrt{am},\sqrt{a\Delta}).
		\]
	\end{prop}
	
	\begin{proof}
		Let $B=\sum_j\ZZ[L_j]$,   $\Lambda=B^\perp$, and let
		$F=L_2+L_3$ be the fibre class of the conic pencil through $L_1$.
		In the original coordinates, choose one component in each of the
		four finite degenerate geometric fibres as follows:
		\begin{equation}\label{eq:alg-fibre-components}
			\begin{aligned}
				C_\pm&:\quad \Delta x=\pm2q_b q_c w,
				& q_b y\mp q_c z&=\Delta w,\\
				Q_\pm&:\quad q_a x=\pm q_m w,
				& 2b\,q_a y&=(\pm q_m+\Delta)z.
			\end{aligned}
		\end{equation}
		In each case, replacing $\Delta$ by $-\Delta$ in the right-hand
		side of the second equation gives the complementary component,
		whose class is $F-C_\pm$ or $F-Q_\pm$, respectively.  Define
		\begin{equation}\label{eq:alg-epsilon-generators}
			\begin{aligned}
				\varepsilon_1&=F-Q_+-Q_-,&
				\varepsilon_2&=F-C_+-C_-,\\
				\varepsilon_3&=C_--C_+,&
				\varepsilon_4&=Q_--Q_+.
			\end{aligned}
		\end{equation}
		Each chosen component meets $L_1$ once and avoids $L_2,L_3$;
		components in distinct fibres are disjoint.  Hence
		$\varepsilon_i\in\Lambda$ and
		$\varepsilon_i\cdot\varepsilon_j=-2\delta_{ij}$.  Moreover,
		\begin{equation}\label{eq:alg-omega}
			\omega:=F-C_+-Q_+
			=\tfrac12(\varepsilon_1+\varepsilon_2+
			\varepsilon_3+\varepsilon_4)
			\in\Lambda.
		\end{equation}
		The primitive boundary lattice $B$ has determinant $4$.
		Put $P=\Pic(\overline X_m)$ and $B^\vee=\Hom(B,\ZZ)$.
		Primitivity of $B$ means that $P/B$ is torsion-free, so a basis
		of $B$ extends to a basis of $P$.  Every homomorphism
		$B\to\ZZ$ extends to $P\to\ZZ$.  Unimodularity of $P$ means
		that its intersection matrix has determinant $\pm1$; equivalently,
		intersection identifies $P$ with $P^\vee=\Hom(P,\ZZ)$.
		The composite $P\simeq P^\vee\to B^\vee$ is therefore surjective.
		It sends $p$ to $(b\mapsto p\cdot b)$ and has kernel
		$B^\perp=\Lambda$.  Hence
		$P/(B\oplus\Lambda)\simeq B^\vee/B$, so
		$[P:B\oplus\Lambda]=|\det B|=4$.  The determinant--index
		formula \cite[Chapter~14, (0.1)--(0.2)]{Huy16} gives
		\[
		|\det B|\,|\det\Lambda|
		=[P:B\oplus\Lambda]^2|\det P|=16,
		\]
		whence $|\det\Lambda|=4$.
		Put $\Lambda_0=\sum_i\ZZ\varepsilon_i$ and
		$\Lambda'=\Lambda_0+\ZZ\omega\subseteq\Lambda$.
		The orthogonal basis $\varepsilon_i$ gives $|\det\Lambda_0|=16$.
		Since $2\omega\in\Lambda_0$ but $\omega\notin\Lambda_0$,
		we have $[\Lambda':\Lambda_0]=2$.  The same determinant--index
		formula yields
		\[
		|\det\Lambda'|=
		\frac{|\det\Lambda_0|}{[\Lambda':\Lambda_0]^2}
		=\frac{16}{4}=4=|\det\Lambda|.
		\]
		Applying it once more to the full-rank inclusion
		$\Lambda'\subseteq\Lambda$ gives
		$[\Lambda:\Lambda']^2=|\det\Lambda'|/|\det\Lambda|=1$.
		Therefore
		\[
		\Lambda=\Lambda'=\sum_i\ZZ\varepsilon_i+\ZZ\omega.
		\]
		Since $\Lambda$ is primitive in the unimodular lattice $P$,
		every homomorphism $\Lambda\to\ZZ$ extends to $P$ and is given
		by intersection with an element of $P$.  The intersection
		map $P\to\Lambda^\vee:=\Hom(\Lambda,\ZZ)$ is therefore surjective.
		Its kernel is $P\cap(B\otimes\QQ)=B$, since $B$ is primitive.
		By \eqref{eq:alg-pic-quotient}, this gives a Galois-equivariant
		isomorphism
		\[
		\Pic(\overline U_m)=P/B\xrightarrow{\sim}\Lambda^\vee,
		\qquad [p]\longmapsto(\lambda\mapsto p\cdot\lambda).
		\]
		For $f\in\Lambda^\vee$, put $z_i=f(\varepsilon_i)\in\ZZ$.
		The relation $2\omega=\sum_i\varepsilon_i$ gives
		$2f(\omega)=\sum_i z_i$, so $\sum_i z_i$ is even.
		Conversely, any integral tuple with even sum defines a unique
		$f\in\Lambda^\vee$ by
		\[
		f(\varepsilon_i)=z_i,\qquad f(\omega)=\tfrac12\sum_i z_i,
		\]
		since these values respect the sole relation among the generators
		of $\Lambda$.  Hence $f\mapsto(z_i)$ identifies
		$\Lambda^\vee$ with the asserted lattice $\mathcal L$.
		
		Let $\sigma_a,\sigma_b,\sigma_c,\sigma_m$ denote the formal
		sign changes of $q_a,q_b,q_c,q_m$, respectively, each fixing the
		other three symbols.  These need not be automorphisms of $K$:
		the actual Galois action is obtained by restricting the sign action
		along
		\[
		\Gal(K/k)\hookrightarrow\{\pm1\}^4,\qquad
		g\longmapsto\left(\frac{g(q_a)}{q_a},\frac{g(q_b)}{q_b},
		\frac{g(q_c)}{q_c},\frac{g(q_m)}{q_m}\right).
		\]
		For example, if $q_a\in k^\times$, the first sign is always $+1$.
		
		We compute the action of two of these formal sign changes, using the
		same notation for divisors and their classes.  The equations give
		\[
		\sigma_a(C_\pm)=C_\pm,\qquad
		\sigma_a(Q_\pm)=F-Q_\mp.
		\]
		It follows that
		\[
		\sigma_a(\varepsilon_1)
		=F-(F-Q_-)-(F-Q_+)=-\varepsilon_1,
		\]
		while $\varepsilon_2,\varepsilon_3,\varepsilon_4$ are fixed.
		Likewise,
		\[
		\sigma_m(C_\pm)=C_\pm,\qquad \sigma_m(Q_\pm)=Q_\mp,
		\]
		so $\sigma_m(\varepsilon_4)=Q_+-Q_-=-\varepsilon_4$,
		while the other three $\varepsilon_i$ are fixed.
		The analogous computations show that $\sigma_b$ and $\sigma_c$
		negate only $\varepsilon_2$ and $\varepsilon_3$, respectively.
		Restricting to the actual Galois group therefore gives exactly the
		characters in \eqref{eq:alg-characters}; the dual action has the
		same characters.  Finally, adding the first character to the other
		three gives $[ab],[ac],[am]$, proving the equality of the splitting
		fields.
	\end{proof}
	
	\begin{lem}\label{lem:alg-dual-basis}
		In the notation of Proposition~\ref{prop:alg-d4},
		\[
		\Lambda=\mathcal L^\vee=\ZZ^4+\ZZ\omega,
		\qquad \omega=(1,1,1,1)/2.
		\]
		An integral basis is
		\[
		A_i=\omega-\varepsilon_i-\varepsilon_4\quad(1\le i\le3),
		\qquad A_4=-\varepsilon_4.
		\]
	\end{lem}
	
	\begin{proof}
		The first assertion was proved above.  The coordinate matrix whose
		columns are $A_1,\ldots,A_4$ relative to the integral basis
		$(\varepsilon_1,\varepsilon_2,\varepsilon_3,\omega)$ has determinant
		$-1$, so the $A_i$ form an integral basis of $\Lambda$.
	\end{proof}
	
	\subsection{Computation of the algebraic part}

	Let $N$ denote the number of nonzero characters in
	\eqref{eq:alg-characters}, and let
	\[
	d=\dim_{\FF_2}\langle\chi_1,\chi_2,\chi_3,\chi_4\rangle.
	\]
	A relation has \emph{odd weight} if it has the form
	$\sum_i\epsilon_i\chi_i=0$ with
	$\sum_i\epsilon_i=1$ in $\FF_2$.
	
	\begin{thm}\label{thm:alg-brauer}
		The algebraic part is given by
		\[
		\Br_1(U_m)/\Br(k)\simeq
		\begin{cases}
			(\ZZ/2)^{N-d},&N<4,\\
			(\ZZ/2)^{5-d},&N=4\text{ and no odd-weight relation exists},\\
			\ZZ/4\oplus(\ZZ/2)^{3-d},&N=4\text{ and an odd-weight relation exists}.
		\end{cases}
		\]
		In the last case, $d\in\{2,3\}$.  If $N=4$, an odd-weight relation
		exists precisely when at least one of
		\[
		abc\Delta,\qquad abm\Delta,\qquad acm\Delta,\qquad bcm\Delta
		\]
		is a square in $k$.  Thus the possible algebraic part is isomorphic to one of 
		$0$, $\ZZ/2$, $(\ZZ/2)^2$, $(\ZZ/2)^3$, $(\ZZ/2)^4$,
		$\ZZ/4$ and $\ZZ/4\oplus\ZZ/2$.
	\end{thm}
	
	\begin{proof}
		By Lemma~\ref{lem:descent} and Proposition~\ref{prop:alg-d4},
		it suffices to compute $H^1(k,\mathcal L)$.  Put
		$\mathcal M=\bigoplus_i\ZZ(\chi_i)$.  By the explicit description
		of $\mathcal L$ in Proposition~\ref{prop:alg-d4}, the coordinate-sum
		map modulo $2$ gives the exact sequence of Galois modules
		\begin{equation}\label{eq:alg-parity-sequence}
			0\longrightarrow\mathcal L\longrightarrow\mathcal M
			\longrightarrow\ZZ/2\longrightarrow0.
		\end{equation}
		The last map sends $(z_i)$ to $\sum_i z_i\bmod2$; it is
		Galois-equivariant with trivial action on $\ZZ/2$, since sign
		changes disappear modulo $2$.
		Here $H^1(k,\mathcal M)=(\FF_2)^N$, and the map
		$\Phi:H^1(k,\mathcal M)\to H^1(k,\ZZ/2)$ sends the $i$th
		basis vector to $\chi_i$.
		If $N<4$, a trivial sign summand makes
		$\mathcal M^{G_k}\to\ZZ/2$ surjective.  Exactness therefore
		identifies $H^1(k,\mathcal L)$ with $\ker\Phi$.
		Since $\operatorname{im}\Phi=\langle\chi_1,\ldots,\chi_4\rangle$
		has dimension $d$, this kernel has dimension $N-d$ over $\FF_2$.
		We obtain $H^1(k,\mathcal L)\simeq(\ZZ/2)^{N-d}$.
		
		If $N=4$, then $\mathcal M^{G_k}=0$, giving
		\begin{equation}\label{eq:alg-extension}
			\begin{gathered}
				0\longrightarrow\ZZ/2\langle\kappa\rangle
				\longrightarrow H^1(k,\mathcal L)\longrightarrow R\longrightarrow0,\\
				R=\ker\bigl((\FF_2)^4\xrightarrow{(\chi_i)}
				H^1(k,\ZZ/2)\bigr).
			\end{gathered}
		\end{equation}
		For $\epsilon=(\epsilon_i)\in R\subset(\FF_2)^4$, choose its
		integer lift $\widetilde\epsilon=(\widetilde\epsilon_i)\in\mathcal M$
		with $\widetilde\epsilon_i\in\{0,1\}$.  Also represent each
		$\chi_i(g)$ by $0$ or $1$.  The action on the $i$th coordinate
		is multiplication by $(-1)^{\chi_i(g)}$, so define
		\[
		c_\epsilon(g):=\frac{g\widetilde\epsilon-\widetilde\epsilon}{2}
		=-(\widetilde\epsilon_i\chi_i(g))_{i=1}^4.
		\]
		Here division by $2$ is applied to an even integer vector.
		Since $\epsilon\in R$, the sum
		$\sum_i\widetilde\epsilon_i\chi_i(g)$ is even for every $g$;
		hence $c_\epsilon(g)\in\mathcal L$.  Moreover, the displayed
		formula gives
		\[
		c_\epsilon(gh)=c_\epsilon(g)+g\,c_\epsilon(h),
		\]
		so $c_\epsilon$ is a cocycle.  For each nonzero $\chi_i$, the
		cocycle $g\mapsto-\chi_i(g)$ represents the generator of
		$H^1(k,\ZZ(\chi_i))\simeq\ZZ/2$.  The image of
		$[c_\epsilon]$ in $H^1(k,\mathcal M)\simeq(\FF_2)^4$ is
		$\epsilon$.
		
		Finally, $2c_\epsilon(g)=g\widetilde\epsilon-\widetilde\epsilon$.
		It follows that $2[c_\epsilon]$ is the connecting image of the parity
		$\sum_i\widetilde\epsilon_i\bmod2$.  Since $\kappa$ is the
		connecting image of $1\in\ZZ/2$, we obtain
		\begin{equation}\label{eq:alg-doubling}
			2[c_\epsilon]=\left(\sum_i\epsilon_i\bmod2\right)\kappa.
		\end{equation}
		The extension is therefore elementary abelian of rank $5-d$
		if all relations have even weight; otherwise it is
		$\ZZ/4\oplus(\ZZ/2)^{3-d}$.
		
		An odd relation among four nonzero characters has weight
		three.  Its three characters span a two-dimensional space, so
		$d\in\{2,3\}$.  Their products give precisely the four stated
		square conditions.
	\end{proof}
	
	\begin{rmk}\label{rmk:alg-generators}
		The proof gives generators as Galois cocycles, including in the
		order-four cases. For each pair of equal square classes among
		$\chi_1,\ldots,\chi_4$, the following quaternion representatives
		adapt the construction of
		\cite[Proposition~3.2 and Theorem~3.4]{CTWX20} to coefficient twists. 
		When $a=b=c=1$, they recover the classes constructed there.
		In the notation of \eqref{eq:alg-characters}, they are
		\[
		\begin{array}{c|c}
			\chi_1=\chi_2 &(z-2\sqrt{ab},a\Delta)\\
			\chi_1=\chi_3 &(y-2\sqrt{ac},a\Delta)\\
			\chi_2=\chi_3 &(x-2\sqrt{bc},b\Delta)\\
			\chi_1=\chi_4 &(x-\sqrt{m/a},a\Delta)\\
			\chi_2=\chi_4 &(y-\sqrt{m/b},b\Delta)\\
			\chi_3=\chi_4 &(z-\sqrt{m/c},c\Delta).
		\end{array}
		\]
		Under the indicated equality, each displayed square root belongs to
		$k$. If the second entry of a quaternion symbol is nonsquare, the two
		components of the corresponding zero divisor are conjugate over the
		associated quadratic extension, so the quaternion is unramified on
		$U_m$. If that entry is a square, the quaternion is trivial.

	\end{rmk}
	
	\section{The full affine Brauer group}\label{sec:complete-affine}
	
	Throughout this section, $\Delta=m-4abc\ne0$ and $m\ne0$. 
	We describe the $n$-torsion of $\Br(U_m)/\Br(k)$ using finite
	Galois modules and recover the full Brauer group by adjoining
	the constant classes.
	
	Retain the square roots $q_a,q_b,q_c,q_m$ chosen in
	Section~\ref{subsec:geometric-configuration}, so that
	\[
	q_a^2=a\Delta,\qquad q_b^2=b\Delta,\qquad
	q_c^2=c\Delta,\qquad q_m^2=m\Delta,
	\]
	and put
	\[
	L=k(q_a,q_b,q_c,q_m),\quad
	M=\frac m{abc},\quad D=\frac\Delta{abc},\quad
	r=\frac{\Delta q_m}{q_aq_bq_c},\quad
	s=\frac{\Delta^2}{q_aq_bq_c},\quad
	\eta=\frac{r-s}{2}.
	\]
	These definitions give $r^2=M$, $s^2=D$, $M-D=4$, and
	\begin{equation}\label{eq:eta-inverse}
		\eta^{-1}=\frac{r+s}{2}.
	\end{equation}
	Denote the four sign characters by
	\[
	\epsilon_i:G_k\longrightarrow\{1,-1\}\qquad(1\le i\le4).
	\]
	They correspond to $a\Delta,b\Delta,c\Delta,m\Delta$, respectively.
	
	\subsection{The boundary lattice and its period}
	
	Let $\mathscr D=D_x+D_y+D_z=X_m\setminus U_m$ be the boundary
	triangle.  Its three components and its three vertices are defined
	over $k$.  Set
	\[
	\Lambda=\{\lambda\in\operatorname{Pic}(\bar X_m):
	\lambda\cdot D_x=\lambda\cdot D_y
	=\lambda\cdot D_z=0\}.
	\]
	The pair $(\bar X_m,\bar{\mathscr D})$ is a Looijenga pair:
	$\bar X_m$ is rational and $\bar{\mathscr D}$ is an anticanonical
	cycle of rational curves.  We recall the algebraic definition of
	its period homomorphism from \cite[Definition~1.5]{GHK15}.
	
	\begin{defn}\label{def:boundary-period}
		The \emph{period homomorphism} of $(\bar X_m,\bar{\mathscr D})$ is
		given by restriction of line bundles:
		\[
		\phi:\Lambda\longrightarrow\operatorname{Pic}^0(\bar{\mathscr D}),
		\qquad [\mathcal L]\longmapsto[\mathcal L|_{\bar{\mathscr D}}].
		\]
		Here $\operatorname{Pic}^0(\bar{\mathscr D})$ consists of line
		bundles having degree zero on each boundary component.
		The intersection conditions defining $\Lambda$ ensure that
		restriction lands in this group.  An orientation of the boundary
		cycle identifies this group with $\kbar^*$, as recalled below.
	\end{defn}
	
	We recall the gluing description of \cite[Lemma~2.1]{GHK15}.
	Orient the boundary cycle as $D_x\to D_y\to D_z\to D_x$, and
	write $P_{ij}=D_i\cap D_j$ for its vertices.  Each component is
	a projective line, so for
	$\mathcal L\in\operatorname{Pic}^0(\bar{\mathscr D})$ we can
	choose nowhere-vanishing sections $\sigma_i$ of
	$\mathcal L|_{D_i}$.  The product
	\[
	\theta(\mathcal L)=
	\frac{\sigma_y(P_{xy})}{\sigma_x(P_{xy})}\,
	\frac{\sigma_z(P_{yz})}{\sigma_y(P_{yz})}\,
	\frac{\sigma_x(P_{zx})}{\sigma_z(P_{zx})}
	\]
	identifies $\operatorname{Pic}^0(\bar{\mathscr D})$ with $\kbar^*$;
	the ratios are taken in the common fibre at each vertex.
	Multiplying $\sigma_i$ by a nonzero constant leaves this product unchanged.
	Conversely, the three gluing scalars, modulo such changes of
	component trivializations, are determined by their product.
	Reversing the orientation inverts $\theta$.
	We henceforth use $\theta$ to regard $\phi$ as a homomorphism
	$\Lambda\to\kbar^*$.
	
	Every point of $D_x\setminus\{P_{xy},P_{xz}\}$ is a smooth
	point of $\bar{\mathscr D}$ and hence defines a
	Cartier divisor on it.  Therefore, any divisor $A=\sum_Pn_PP$ supported
	in $D_x\setminus\{P_{xy},P_{xz}\}$ is a Cartier divisor on the
	boundary.  The notation
	$\mathcal O_{\bar{\mathscr D}}(A)$ denotes its associated line
	bundle, which is canonically trivial away from $\operatorname{Supp}(A)$,
	and in particular on $D_y,D_z$.
	
	\begin{lem}\label{lem:boundary-period}
		Use the coordinate $u=Y/Z$ on $D_x$, so that $P_{xy}$ and $P_{xz}$ have
		coordinates $0$ and $\infty$.  Let $A=\sum_Pn_PP$ be a
		degree-zero Cartier divisor on $\bar{\mathscr D}$ supported in
		$D_x\setminus\{P_{xy},P_{xz}\}$.  For any rational function
		$f$ on $D_x$ with $\operatorname{div}_{D_x}(f)=A|_{D_x}$, one has
		\[
		\theta\bigl(\mathcal O_{\bar{\mathscr D}}(A)\bigr)
		=\frac{f(0)}{f(\infty)}=\prod_Pu(P)^{n_P}.
		\]
	\end{lem}
	
    \begin{proof}
    	Put $\mathcal L=\mathcal O_{\bar{\mathscr D}}(A)$ and use the
    	component sections $\sigma_i$ in the gluing formula above.
    	Since $\operatorname{div}_{D_x}(f)=A|_{D_x}$, the divisor-sheaf
    	description gives
    	$\mathcal L|_{D_x}=f^{-1}\mathcal O_{D_x}$.  On $D_y,D_z$,
    	the line bundle is canonically trivial.  We may therefore choose
    	\[
    	\sigma_x=f^{-1},\qquad \sigma_y=1,\qquad \sigma_z=1.
    	\]
    	These are nowhere-vanishing sections of the respective line
    	bundles.  At the vertices, $A$ has no support, so the sections
    	can be compared using the common local trivialization $1$.  We obtain
    	\[
    	\frac{\sigma_y(P_{xy})}{\sigma_x(P_{xy})}=f(0),\qquad
    	\frac{\sigma_z(P_{yz})}{\sigma_y(P_{yz})}=1,\qquad
    	\frac{\sigma_x(P_{zx})}{\sigma_z(P_{zx})}=f(\infty)^{-1}.
    	\]
    	The formula defining $\theta$ gives
    	$\theta(\mathcal L)=f(0)/f(\infty)$.
    	Taking
    	\[
    	f(u)=\prod_P\bigl(u-u(P)\bigr)^{n_P}
    	\]
    	gives $f(\infty)=1$ and $f(0)=\prod_Pu(P)^{n_P}$, because
    	$\sum_Pn_P=0$.  Any other choice differs by a nonzero constant,
    	which preserves the ratio.
    \end{proof}
	
	\begin{prop}\label{prop:affine-period}
		As a Galois lattice, $\Lambda$ has the following coordinate model:
		\[
		\Lambda=\mathbb Z e_1+\mathbb Z e_2+\mathbb Z e_3
		+\mathbb Z e_4+\mathbb Z h,
		\qquad 2h=e_1+e_2+e_3+e_4.
		\]
		Here the $e_i$'s are linearly independent in $\Lambda\otimes\mathbb Q$ with Galois action 
		$g(e_i)=\epsilon_i(g)e_i$.  The period is
		\begin{equation}\label{eq:period-weight-basis}
			\phi(e_1)=\phi(e_2)=\phi(e_3)=-1,\qquad
			\phi(e_4)=-\eta^2,\qquad \phi(h)=\eta.
		\end{equation}
		Moreover,
		\[
		\operatorname{Pic}(\bar U_m)\simeq\Lambda^\vee
		:=\operatorname{Hom}(\Lambda,\mathbb Z).
		\]
	\end{prop}
	
	\begin{proof}
		Take $e_i=\varepsilon_i$ and $h=\omega$ from the proof of
		Proposition~\ref{prop:alg-d4}, which establishes the lattice
		structure, the Galois action and the duality assertion.
		It remains to apply Lemma~\ref{lem:boundary-period} to the
		boundary intersections of their defining divisors.
		
		We work in the normalized coordinates
		\eqref{eq:alg-coordinate-change}, retaining the notation for the
		transported components.
		Represent the fibre class $F=D_y+D_z$ by the smooth conic
		\[
		F_0:\quad X=0,\qquad Y^2+Z^2=Mw^2.
		\]
		Indeed, the plane sections $X=0$ and $w=0$ are $D_x+F_0$ and
		$D_x+D_y+D_z$ respectively, so $[F_0]=F$.
		Its two intersection coordinates on $D_x$ satisfy $u^2+1=0$,
		so their product is $1$.  In these coordinates, the component
		equations \eqref{eq:alg-fibre-components} become
		\[
		\begin{aligned}
			C_\pm&:\quad X=\pm2w,&Y\mp Z&=sw,\\
			Q_\pm&:\quad X=\pm rw,&Y&=\frac{\pm r+s}{2}Z.
		\end{aligned}
		\]
		Setting $w=0$ gives the $u$-coordinates
		$1,-1,(r+s)/2=\eta^{-1},(-r+s)/2=-\eta$, respectively.
		All these curves meet the boundary only on $D_x$, transversely and
		away from the vertices.  We obtain the following products of their
		intersection coordinates, counted with multiplicities:
		\[
		\begin{array}{c|ccccc}
			E&F_0&C_+&C_-&Q_+&Q_-\\ \hline
			\displaystyle\prod_{P\in E\cap D_x}u(P)&1&1&-1&\eta^{-1}&-\eta
		\end{array}
		\]
		For any of the surface divisors $E$ defining
		$\varepsilon_i,\omega$, with $F$ represented by $F_0$, its support
		is disjoint from $D_y\cup D_z$.  Hence
		$A:=E|_{\bar{\mathscr D}}$ is a Cartier divisor supported in
		$D_x\setminus\{P_{xy},P_{xz}\}$, and
		\[
		\mathcal O_{\bar X_m}(E)|_{\bar{\mathscr D}}
		\simeq\mathcal O_{\bar{\mathscr D}}(A).
		\]
		Since $[E]\in\Lambda$, this divisor $A$ has degree zero.
		Lemma~\ref{lem:boundary-period} therefore applies and gives, for example,
		\[
		\phi(\omega)=\frac{1}{1\cdot\eta^{-1}}=\eta,
		\qquad
		\phi(\varepsilon_4)=\frac{-\eta}{\eta^{-1}}=-\eta^2.
		\]
		The same formula applied to $\varepsilon_1=F-Q_+-Q_-$,
		$\varepsilon_2=F-C_+-C_-$ and $\varepsilon_3=C_--C_+$ gives
		$\phi(\varepsilon_1)=\phi(\varepsilon_2)=
		\phi(\varepsilon_3)=-1$.  With the generators defined by
		\eqref{eq:alg-epsilon-generators} and \eqref{eq:alg-omega},
		and the identifications $e_i=\varepsilon_i$, $h=\omega$, these
		are precisely the period values in \eqref{eq:period-weight-basis}.
	\end{proof}
	
	The period homomorphism is Galois-equivariant: restriction to the
	boundary and the gluing identification
	$\Pic^0(\bar{\mathscr D})\simeq\kbar^*$ are defined over $k$.
	
	\subsection{A finite module computing the full group}
	
	For a positive integer $n$, we define the finite Galois module
	\begin{equation}\label{eq:definition-Tn}
		T_n=
		\frac{\{(\lambda,t)\in\Lambda\times\kbar^*:
			t^n=\phi(\lambda)\}}
		{\{(n\lambda,\phi(\lambda)):\lambda\in\Lambda\}}.
	\end{equation}
	The group law is addition in the first coordinate and
	multiplication in the second.  It fits into an exact sequence
	\begin{equation}\label{eq:Tn-extension}
		0\longrightarrow\mu_n\longrightarrow T_n
		\longrightarrow\Lambda/n\Lambda\longrightarrow0.
	\end{equation}
	Every element is killed by $n$, and $T_n$ is free of rank five
	over $\mathbb Z/n$.  Indeed, multiplication by $n$ sends
	$(\lambda,t)$ to the defining relation $(n\lambda,\phi(\lambda))$;
	the quotient $\Lambda/n\Lambda$ is free of rank four, so
	\eqref{eq:Tn-extension} splits as a sequence of $\mathbb Z/n$-modules.
	Its Galois extension class is part of the calculation below.
	
	\begin{thm}\label{thm:complete-affine-brauer}
		The inclusion $\operatorname{Br}(k)\to\operatorname{Br}(U_m)$
		is split by evaluation at a $k$-point of $U_m$.  For every
		$n\ge1$, there is an isomorphism
		\begin{equation}\label{eq:full-brauer-finite}
			\left(\operatorname{Br}(U_m)/\operatorname{Br}(k)\right)[n]
			\simeq
			\frac{\operatorname{Hom}_{G_k}(T_n,\mathbb Z/n)}
			{ (\Lambda^\vee)^{G_k}/n(\Lambda^\vee)^{G_k}}.
		\end{equation}
		The denominator maps to the numerator by
		$\ell(\lambda,t)=\ell(\lambda)\bmod n$.  Taking the direct limit
		with respect to the maps induced by the Kummer sequences computes the
		full group modulo constants.
	\end{thm}
	
	For the direct limit appearing in Theorem~\ref{thm:complete-affine-brauer},
	if $n\mid N$ and $d=N/n$, the reduction map
	$T_N\to T_n$ sends $(\lambda,t)$ to $(\lambda,t^d)$.
	The transition on the numerator sends a functional $F$ to
	$dF\circ\mathrm{red}$ with values in $\ZZ/N$; on the denominator
	it sends $\ell\bmod n$ to $d\ell\bmod N$.  These maps specify the
	direct limit, so evaluation at a rational point gives
	\[
	\Br(U_m)\simeq\Br(k)\oplus
	\varinjlim_n
	\frac{\Hom_{G_k}(T_n,\ZZ/n)}
	{ (\Lambda^\vee)^{G_k}/n(\Lambda^\vee)^{G_k}}.
	\]
	
	\medskip
	We prove Theorem~\ref{thm:complete-affine-brauer} using three lemmas.
	Fix $n\ge1$, put $R=\ZZ/n$, and write
	$i:\bar{\mathscr D}\hookrightarrow\bar X_m$ and
	$j:\bar U_m\hookrightarrow\bar X_m$ for the inclusions.
	All geometric cohomology groups carry their natural $G_k$-action.
	
	\begin{lem}\label{lem:boundary-cohomology}
		The multidegree map gives an exact sequence
		\begin{equation}\label{eq:boundary-picard-degree}
			0\longrightarrow \kbar^*\longrightarrow
			\operatorname{Pic}(\bar{\mathscr D})
			\xrightarrow{\deg}\ZZ^3\longrightarrow0.
		\end{equation}
		Its kernel is $\operatorname{Pic}^0(\bar{\mathscr D})$,
		and $\operatorname{Br}(\bar{\mathscr D})=0$.  Moreover,
		\begin{equation}\label{eq:boundary-finite-cohomology}
			H^1(\bar{\mathscr D},\mu_n)=\kbar^*[n]=\mu_n,\qquad
			H^2(\bar{\mathscr D},\mu_n)=R^3.
		\end{equation}
		The relative cohomology sequence of the pair
		$(\bar X_m,\bar{\mathscr D})$ induces an exact sequence
		\begin{equation}\label{eq:relative-multidegree-extension}
			0\longrightarrow\mu_n\longrightarrow H^2_c(\bar U_m,\mu_n)
			\longrightarrow\Lambda/n\Lambda\longrightarrow0.
		\end{equation}
	\end{lem}
	
	\begin{proof}
		Let $\nu:\widetilde{\mathscr D}\to\bar{\mathscr D}$ be the normalization
		and let $p$ run through the three nodes.  The normalization sequence is
		\begin{equation}\label{eq:boundary-normalization-units}
			1\longrightarrow\mathbf G_{m,\bar{\mathscr D}}
			\longrightarrow\nu_*\mathbf G_{m,\widetilde{\mathscr D}}
			\longrightarrow\bigoplus_p(i_p)_*\mathbf G_{m,p}
			\longrightarrow1.
		\end{equation}
		The relevant higher direct images of $\nu$ vanish on strict local
		branches.  On invertible regular functions, the last map is
		$(c_x,c_y,c_z)\mapsto(c_y/c_x,c_z/c_y,c_x/c_z)$, with cokernel $\kbar^*$.
		Since $\operatorname{Pic}(\widetilde{\mathscr D})=\ZZ^3$ and
		$\operatorname{Pic}(p)=\operatorname{Br}(\widetilde{\mathscr D})=0$,
		its cohomology gives \eqref{eq:boundary-picard-degree} and
		$\operatorname{Br}(\bar{\mathscr D})=0$.
		The kernel of multidegree is the connected torus $\mathbf G_m$
		of cyclic gluing parameters; since multidegree is locally constant,
		this kernel equals $\operatorname{Pic}^0(\bar{\mathscr D})$.
		The Kummer sequence and the divisibility of $\kbar^*$ now give
		\eqref{eq:boundary-finite-cohomology}.
		
		Properness of $\bar X_m$ and the open--closed sequence
		$0\to j_!\mu_n\to\mu_{n,\bar X_m}\to i_*\mu_{n,\bar{\mathscr D}}\to0$
		give
		\begin{equation}\label{eq:relative-etale-localization}
			\begin{aligned}
				H^1(\bar X_m,\mu_n)&\longrightarrow
				H^1(\bar{\mathscr D},\mu_n)\longrightarrow
				H^2_c(\bar U_m,\mu_n)\\
				&\longrightarrow H^2(\bar X_m,\mu_n)
				\longrightarrow H^2(\bar{\mathscr D},\mu_n).
			\end{aligned}
		\end{equation}
		Since $\bar X_m$ is rational,
		$H^1(\bar X_m,\mu_n)=0$ and
		$H^2(\bar X_m,\mu_n)=\operatorname{Pic}(\bar X_m)/n$.
		The last map is multidegree modulo $n$.  By
		Proposition~\ref{prop:alg-d4}, the integral multidegree map
		$\Pic(\bar X_m)\to\ZZ^3$ is surjective with kernel $\Lambda$.
		Since $\ZZ^3$ is free, reduction modulo $n$ gives an exact sequence
		\[
		0\longrightarrow\Lambda/n\Lambda
		\longrightarrow\Pic(\bar X_m)/n
		\longrightarrow(\ZZ/n)^3\longrightarrow0.
		\]
		The last map in \eqref{eq:relative-etale-localization}
		therefore has kernel $\Lambda/n\Lambda$, proving
		\eqref{eq:relative-multidegree-extension}.
	\end{proof}
	
	\begin{lem}\label{lem:relative-Kummer}
		Let
		\[
		\mathcal F=\ker\bigl(\mathbf G_{m,\bar X_m}
		\longrightarrow i_*\mathbf G_{m,\bar{\mathscr D}}\bigr).
		\]
		In degrees one and two, we have
		\[
		H^1(\bar X_m,\mathcal F)=\ker(\phi:\Lambda\to\kbar^*),
		\qquad H^2(\bar X_m,\mathcal F)=\kbar^*/\phi(\Lambda).
		\]
		For every $n\ge1$, the relative Kummer sequence identifies
		\eqref{eq:relative-multidegree-extension} with \eqref{eq:Tn-extension}
		by a $G_k$-equivariant isomorphism
		\begin{equation}\label{eq:compact-realization}
			H^2_c(\bar U_m,\mu_n)\simeq T_n.
		\end{equation}
	\end{lem}
	
	\begin{proof}
		Write $X=\bar X_m$, $D=\bar{\mathscr D}$ and $U=\bar U_m$.
		The complexes below are
		complexes of discrete $G_k$-modules, with cohomological shifts
		$H^q(K[r])=H^{q+r}(K)$.  We write $\tau_{\le r}$ for the canonical
		truncation preserving cohomology in degrees at most $r$.
		
		\smallskip\noindent
		\emph{Relative invertible regular functions and line bundles.}
		Invertible regular functions lift locally across a closed immersion, so
		\begin{equation}\label{eq:relative-F-units}
			1\longrightarrow\mathcal F\longrightarrow\mathbf G_{m,X}
			\longrightarrow i_*\mathbf G_{m,D}\longrightarrow1
		\end{equation}
		is exact.  A class in $H^1(X,\mathcal F)$ is represented by a line
		bundle $\mathcal L$ on $X$ together with a trivialization
		$\alpha:\mathcal L|_D\xrightarrow{\sim}\mathcal O_D$.
		Indeed, local frames compatible with $\alpha$ have transition
		functions restricting to $1$ on $D$, and conversely such transition
		functions glue a pair $(\mathcal L,\alpha)$.
		This is the relative Picard group
		$\operatorname{Pic}(X,D)$; see
		\cite[Section~2.1 and Proposition~2.5.1]{BVS01}.
		
		The groups of invertible regular functions on $X$ and $D$ are both
		$\kbar^*$, and restriction between them is the identity.  The cohomology sequence of
		\eqref{eq:relative-F-units} therefore gives
		\begin{equation}\label{eq:relative-Pic-kernel}
			H^1(X,\mathcal F)=\operatorname{Pic}(X,D)
			=\ker\bigl[\operatorname{Pic}(X)\to\operatorname{Pic}(D)\bigr]
			=\ker(\phi:\Lambda\to\kbar^*).
		\end{equation}
		A line bundle trivial on $D$ has multidegree zero there,
		so its class lies in $\Lambda$, and its restriction is trivial
		exactly when $\phi(\lambda)=1$.
		Since $H^2(X,\mathbf G_m)=0$, the same sequence identifies
		$H^2(X,\mathcal F)$ with the cokernel of restriction on Picard groups.
		The integral multidegree map $\operatorname{Pic}(X)\to\mathbb Z^3$
		is surjective. Thus every line bundle on $D$ is
		represented, modulo restriction from $X$, by an element of
		$\operatorname{Pic}^0(D)=\kbar^*$; the remaining relations are
		precisely $\phi(\Lambda)$.  Hence
		\[
		H^2(X,\mathcal F)=\kbar^*/\phi(\Lambda).
		\]
		
		\smallskip\noindent
		\emph{The relative Kummer exact sequence.}
		For $n\ge1$, multiplication by $n$ on $\mathcal F$ is surjective
		as a map of \'{e}tale sheaves.  On $U$ this is ordinary Kummer.
		Near $D$, if $f|_D=1$, the equation $T^n=f$ is \'{e}tale and its
		restriction to $D$ has the branch $T=1$.  Choosing this branch
		gives an \'{e}tale-local root still restricting to $1$ on $D$,
		including at the nodes.  The kernel is $\mu_n$ on $U$ and has
		zero stalks on $D$: an $n$-th root of unity whose restriction
		to $D$ is $1$ is itself $1$ in the strict local neighbourhood.
		We therefore obtain
		\begin{equation}\label{eq:relative-F-Kummer}
			0\longrightarrow j_!\mu_n\longrightarrow\mathcal F
			\xrightarrow{(\cdot)^n}\mathcal F\longrightarrow0.
		\end{equation}
		This is the relative Kummer sheaf sequence stated at the beginning
		of \cite[Section~2.5, p.~33]{BVS01}.
		Proposition~2.5.2 of \cite{BVS01} describes its degree-one cohomology.
		Now we give the derived calculation needed in degree two.
		
		\smallskip\noindent
		\emph{A model for the relative derived global sections.}
		Put $C=R\Gamma(X,\mathcal F)$ and
		$K_Y=R\Gamma(Y,\mathbf G_m)$ for $Y=X,D$.
		By \eqref{eq:relative-F-units},
		$C\simeq\operatorname{Cone}(K_X\to K_D)[-1]$.
		Evaluation at a common rational boundary vertex splits off
		$(\kbar^*)[0]$ from both $K_X$ and $K_D$, compatibly with
		restriction and the $G_k$-action.  Let $\widetilde K_Y$ denote
		the complementary summands.  Restriction is the identity on
		constant invertible regular functions, whose cone is acyclic, so
		\[
		C\simeq\operatorname{Cone}(\widetilde K_X\to\widetilde K_D)[-1].
		\]
		The reduced complexes have $H^0=H^2=0$ and
		$H^1=\operatorname{Pic}(Y)$, using the vanishing of both geometric
		Brauer groups.  Hence
		$\tau_{\le2}\widetilde K_Y\simeq\operatorname{Pic}(Y)[-1]$.
		Taking the relative fibre and truncating in degrees at most two gives
		\begin{equation}\label{eq:relative-Pic-complex}
			\tau_{\le2}C\simeq
			[\operatorname{Pic}(X)\xrightarrow{\mathrm{res}}
			\operatorname{Pic}(D)]
			\quad\text{in degrees $1,2$}.
		\end{equation}
		The subcomplex
		\[
		C_\phi=[\Lambda\xrightarrow{\phi}\kbar^*]
		\quad\text{in degrees $1,2$}
		\]
		embeds into this Picard complex.  By the surjectivity of the integral multidegree map
		and \eqref{eq:boundary-picard-degree}, the quotient complex is
		$[\mathbb Z^3\xrightarrow{\mathrm{id}}\mathbb Z^3]$, which is acyclic.
		We therefore have $\tau_{\le2}C\simeq C_\phi$ equivariantly, with the
		differential $\phi$ retained.
		
		\smallskip\noindent
		\emph{Applying Kummer and computing degree two.}
		Apply $R\Gamma(X,-)$ to \eqref{eq:relative-F-Kummer}.
		The resulting distinguished triangle is
		\[
		R\Gamma(X,j_!\mu_n)\longrightarrow C
		\xrightarrow{n}C\longrightarrow R\Gamma(X,j_!\mu_n)[1].
		\]
		Therefore
		\[
		R\Gamma(X,j_!\mu_n)
		\simeq\operatorname{Cone}(n:C\to C)[-1].
		\]
		Properness of $X$ identifies the left side with
		$R\Gamma_c(U,\mu_n)$.
		The long exact cohomology sequence shows that we may replace
		$C$ by $\tau_{\le2}C\simeq C_\phi$ when computing the
		degree-two cohomology of this shifted cone.  Hence
		\[
		H^2_c(U,\mu_n)\simeq
		H^2\bigl(\operatorname{Cone}(n:C_\phi\to C_\phi)[-1]\bigr).
		\]
		With compatible cone signs, the explicit complex on the right is
		\begin{equation}\label{eq:relative-Kummer-three-term}
			\Lambda\xrightarrow{d^1}\Lambda\oplus\kbar^*
			\xrightarrow{d^2}\kbar^*
			\quad\text{in degrees $1,2,3$},
		\end{equation}
		where
		\[
		d^1(\nu)=(n\nu,\phi(\nu)),\qquad
		d^2(\lambda,t)=\phi(\lambda)t^{-n}.
		\]
		The composition is trivial because $\phi(n\nu)=\phi(\nu)^n$.
		Taking the kernel of $d^2$ modulo the image of $d^1$ gives
		\[
		H^2_c(U,\mu_n)\simeq
		\frac{\{(\lambda,t)\in\Lambda\times\kbar^*:
			t^n=\phi(\lambda)\}}
		{\{(n\nu,\phi(\nu)):\nu\in\Lambda\}}
		=T_n.
		\]
		Here $t$ represents a boundary line bundle $\mathcal N$ of
		multidegree zero with
		$\mathcal N^{\otimes n}\simeq\mathcal L|_D$, where
		$[\mathcal L]=\lambda$.  
		
		Every step is equivariant, and on representatives the action is
		$g[(\lambda,t)]=[(g\lambda,g(t))]$.
		The maps from the boundary sequence become
		\[
		\mu_n\longrightarrow T_n,\quad\zeta\longmapsto[(0,\zeta)],
		\qquad
		T_n\longrightarrow\Lambda/n\Lambda,\quad
		[(\lambda,t)]\longmapsto\lambda\bmod n\Lambda.
		\]
		They identify \eqref{eq:relative-multidegree-extension} with
		\eqref{eq:Tn-extension}, including its Galois extension class,
		and prove \eqref{eq:compact-realization}.
	\end{proof}
	
	\begin{lem}\label{lem:duality-Chern}
		There is a $G_k$-equivariant isomorphism
		\begin{equation}\label{eq:ordinary-realization}
			H^2(\bar U_m,\mu_n)\simeq
			\operatorname{Hom}(T_n,\ZZ/n).
		\end{equation}
		Under this isomorphism and the identification
		$\operatorname{Pic}(\bar U_m)=\Lambda^\vee$,
		the first Chern class of $\ell$ is the functional
		\[
		[(\lambda,t)]\longmapsto\ell(\lambda)\bmod n.
		\]
		Passing to the quotient gives an induced equivariant isomorphism
		\[
		\Br(\bar U_m)[n]\simeq\Hom(\mu_n,\ZZ/n)
		\]
		under which the Kummer projection is restriction of functionals
		from $T_n$ to $\mu_n$.
	\end{lem}
	
	\begin{proof}
		Poincar\'e duality \cite{Milne} gives a perfect equivariant pairing
		\[
		H^2(\bar U_m,\mu_n)\times H^2_c(\bar U_m,\mu_n)
		\longrightarrow H^4_c(\bar U_m,\mu_n^{\otimes2})
		\xrightarrow{\operatorname{tr}}\ZZ/n.
		\]
		Here $\ZZ/n$ has trivial Galois action; the two coefficient
		twists have already combined to give the surface trace twist.
		Together with Lemma~\ref{lem:relative-Kummer}, this proves
		\eqref{eq:ordinary-realization}.
		Dualizing \eqref{eq:Tn-extension} gives
		\[
		0\longrightarrow\Lambda^\vee/n\Lambda^\vee
		\longrightarrow\operatorname{Hom}(T_n,\ZZ/n)
		\xrightarrow{\mathrm{res}}\Hom(\mu_n,\ZZ/n)
		\longrightarrow0.
		\]
		To identify the left map, lift a line bundle on $\bar U_m$ to a
		divisor class $\widetilde\ell\in\operatorname{Pic}(\bar X_m)$.
		Its image $\ell\in\Lambda^\vee$ is the functional
		$\lambda\mapsto\widetilde\ell\cdot\lambda$; changing the lift
		by boundary classes does not change this functional.
		The map forgetting support sends $[(\lambda,t)]$ to
		$\lambda\bmod n\in H^2(\bar X_m,\mu_n)$.
		Compatibility of cup products, cycle classes and the trace therefore
		shows that its pairing with the first Chern class of $\ell$ is
		$\widetilde\ell\cdot\lambda=\ell(\lambda)\bmod n$.
		The left map is therefore the first Chern class map.
		
		To identify the right-hand map as well, compare this dual sequence
		with the geometric Kummer sequence
		\[
		0\longrightarrow\Pic(\bar U_m)/n
		\xrightarrow{c_1}H^2(\bar U_m,\mu_n)
		\longrightarrow\Br(\bar U_m)[n]\longrightarrow0.
		\]
		The middle terms are identified by duality, and the left terms
		by $\Pic(\bar U_m)=\Lambda^\vee$.
		The Chern class calculation above identifies the two left
		injections.  Passing to their cokernels therefore gives the
		right-hand isomorphism and the commutative diagram
		\begin{equation}\label{eq:brauer-restriction-duality}
			\begin{CD}
				H^2(\bar U_m,\mu_n) @>>> \Br(\bar U_m)[n]\\
				@V{\sim}VV @VV{\sim}V\\
				\Hom(T_n,\ZZ/n) @>{F\mapsto F|_{\mu_n}}>>
				\Hom(\mu_n,\ZZ/n).
			\end{CD}
		\end{equation}
		In particular, the bottom arrow is restriction because it is
		the dual of the inclusion $\mu_n\hookrightarrow T_n$.
		
		Let $F$ correspond by duality to the geometric restriction of
		a Kummer lift of $\beta\in\Br(U_m)[n]$.
		By naturality under base change,
		diagram~\eqref{eq:brauer-restriction-duality} identifies the
		geometric restriction of $\beta$ with $F|_{\mu_n}$.
		This is independent of the lift, since Picard functionals
		vanish on $\mu_n$.
	\end{proof}
	
	\begin{proof}[Proof of Theorem~\ref{thm:complete-affine-brauer}]
		Choose $P_0\in U_m(k)$ as in Lemma~\ref{lem:smooth-points}, and put
		\[
		C_n=\ker\bigl[P_0^*:H^2(U_m,\mu_n)\longrightarrow H^2(k,\mu_n)\bigr].
		\]
		Evaluation splits pullback from $k$, so Lemma~\ref{lem:descent}
		and Lemma~\ref{lem:duality-Chern} give
		\[
		C_n\simeq H^2(\bar U_m,\mu_n)^{G_k}
		\simeq\operatorname{Hom}_{G_k}(T_n,\ZZ/n).
		\]
		Moreover, Lemma~\ref{lem:descent} gives
		\begin{equation}\label{eq:Pic-affine-descent}
			\operatorname{Pic}(U_m)=(\Lambda^\vee)^{G_k}.
		\end{equation}
		First Chern classes lie in $C_n$, because
		$\operatorname{Pic}(\operatorname{Spec}k)=0$.
		By Lemma~\ref{lem:duality-Chern}, the Chern class of $\ell$
		corresponds to $(\lambda,t)\mapsto\ell(\lambda)\bmod n$.
		
		Compatibility of the Kummer sequence with evaluation at $P_0$ therefore gives
		\[
		0\longrightarrow\operatorname{Pic}(U_m)/n
		\longrightarrow C_n
		\longrightarrow
		\ker\bigl[\operatorname{Br}(U_m)\xrightarrow{P_0^*}
		\operatorname{Br}(k)\bigr][n]
		\longrightarrow0.
		\]
		Evaluation splits the entire constant Brauer subgroup, so its
		kernel is isomorphic to
		$\operatorname{Br}(U_m)/\operatorname{Br}(k)$.
		Substituting \eqref{eq:Pic-affine-descent} proves
		\eqref{eq:full-brauer-finite}.  The transition maps follow from
		the functoriality of the relative Kummer sequence and the duality pairing.
		
	\end{proof}
	
	\subsection{An explicit formula for the transcendental part}
	
	We determine the transcendental part using the connecting
	homomorphism in Galois cohomology.  Retain $M=m/(abc)$ and
	$D=\Delta/(abc)$, and recall
	\[
	r=\frac{\Delta q_m}{q_aq_bq_c},\qquad
	s=\frac{\Delta^2}{q_aq_bq_c},\qquad
	\eta=\frac{r-s}{2}.
	\]
	These definitions give $r^2=M$, $s^2=D$ and $rs=q_m/(abc)$.  Retain also
	\[
	L=k(\sqrt{a\Delta},\sqrt{b\Delta},\sqrt{c\Delta},
	\sqrt{m\Delta}).
	\]
	Put
	\[
	K_1=k(\sqrt{a\Delta}),\qquad
	K_2=k(\sqrt{b\Delta}),\qquad
	K_3=k(\sqrt{c\Delta}).
	\]
	Here and below, $\mu_j\subset F$ means that $F$ contains all
	$j$th roots of unity in the fixed algebraic closure.
	
	We first recall the direct-limit form of \eqref{eq:Tn-extension}.
	For $N=dn$, the injections
	\[
	T_n\longrightarrow T_N,\qquad
	[(\lambda,t)]\longmapsto[(d\lambda,t)]
	\]
	induce the usual inclusions on $\mu_n$ and multiplication by $d$
	on $\Lambda/n\Lambda$.  Setting
	$\mathcal T=\varinjlim_n T_n$, we obtain the direct-limit sequence
	\begin{equation}\label{eq:Tn-direct-limit}
		0\longrightarrow\QQ/\ZZ(1)\longrightarrow\mathcal T
		\longrightarrow\Lambda\otimes\QQ/\ZZ\longrightarrow0.
	\end{equation}
	Here $\QQ/\ZZ(1)=\bigcup_n\mu_n$ and $T_n=\mathcal T[n]$,
	so taking $n$-torsion recovers \eqref{eq:Tn-extension}.
	
	For the dual sequence, dualize \eqref{eq:Tn-extension} at each
	finite level and then take direct limits.  Put
	\[
	V_n=\Hom(T_n,\ZZ/n),\qquad
	\mathcal V=\varinjlim_n V_n,\qquad
	\mathcal A=\Lambda^\vee\otimes\QQ/\ZZ.
	\]
	Lemma~\ref{lem:duality-Chern} therefore gives
	\begin{equation}\label{eq:Tn-dual-direct-limit}
		0\longrightarrow\mathcal A\longrightarrow\mathcal V
		\xrightarrow{\mathrm{res}}\QQ/\ZZ(-1)\longrightarrow0.
	\end{equation}
	Here $V_n=\mathcal V[n]$, and $\mathrm{res}$ is induced by
	restriction to $\mu_n$.  In this sense, $\mathcal V$ is the direct-limit
	version of the numerator.
	The injective transition maps give
	$\mathcal V^{G_k}=\varinjlim_n\Hom_{G_k}(T_n,\ZZ/n)$.
	Taking invariants in \eqref{eq:Tn-dual-direct-limit} and using
	Theorem~\ref{thm:complete-affine-brauer} yields
	\begin{equation}\label{eq:trans-obstruction-kernel}
		I:=\frac{\Br(U_m)}{\Br_1(U_m)}
		\simeq\ker\!\left[
		(\QQ/\ZZ(-1))^{G_k}\xrightarrow{\partial}
		H^1(k,\mathcal A)\right].
	\end{equation}
	For an invariant character $t$ and any lift $f\in\mathcal V$,
	\[
	\partial(t)=[g\longmapsto gf-f],\qquad
	(gf)(v)=f(g^{-1}v).
	\]
	This class is independent of the lift: replacing $f$ by $f+a$,
	with $a\in\mathcal A$, adds the coboundary $ga-a$.
	Using the direct limit allows us to consider a Brauer lift whose
	order exceeds that of its geometric image.
	The theorem below computes this kernel.
	
	\begin{thm}\label{trans-thm:explicit}
		Let $k$ be a field of characteristic zero with $a,b,c\in k^*$, and suppose
		$m\ne0,4abc$.  Define
		\begin{equation}\label{trans-eq:index}
			\mathcal I_{a,b,c,m}=\left\{n\ge1:
			\begin{array}{l}
				\mu_n\subset k,\\
				\mu_{2n}\subset K_i\quad\text{for }i=1,2,3,\\
				\eta\in L^{\times n}
			\end{array}\right\}.
		\end{equation}
		Then
		\begin{equation}\label{trans-eq:quotient}
			\frac{\Br(U_m)}{\Br_1(U_m)}
			\simeq\varinjlim_{n\in\mathcal I_{a,b,c,m}}\ZZ/n,
		\end{equation}
		where, for $n\mid n'$, the transition map sends
		$1\bmod n$ to $n'/n\bmod n'$.  The set
		$\mathcal I_{a,b,c,m}$ is closed under divisors and least common
		multiples, and is independent of the choices of square roots.
	\end{thm}
	
	\begin{proof}
		An invariant character of order $n$ requires $\mu_n\subset k$.
		Under this assumption, choose a primitive $n$th root $\zeta_n$
		and write $t_n=1/n$ for the character $t_n(\zeta_n^j)=j/n$.
		Choose a primitive $2n$th root $\xi_{2n}$, an $n$th root
		$\sqrt[n]{\eta}$, and a lift
		$f\in\Hom(T_n,\tfrac1n\ZZ/\ZZ)$ of $t_n$.
		Put $c(g)=gf-f$, so $[c]=\partial(t_n)$.
		
		Suppose $\partial(t_n)=0$, so $c(g)=ga-a$ for some
		$a\in\mathcal A$.
		For $i=1,2,3$ and $g\in G_{K_i}$, the equality $g(e_i)=e_i$
		gives $c(g)(e_i)=a(g^{-1}e_i)-a(e_i)=0$.
		Since $\phi(e_i)=-1$ by \eqref{eq:period-weight-basis}, evaluation
		using the lift $[(e_i,\xi_{2n})]$ gives
		\begin{equation}\label{eq:trans-coordinate-cocycle}
			\begin{aligned}
				c(g)(e_i)
				&=f([(e_i,g^{-1}(\xi_{2n}))])-f([(e_i,\xi_{2n})])\\
				&=f([(0,g^{-1}(\xi_{2n})/\xi_{2n})])\\
				&=t_n\!\left(\frac{g^{-1}(\xi_{2n})}{\xi_{2n}}\right).
			\end{aligned}
		\end{equation}
		The ratio belongs to $\mu_n$ because $\xi_{2n}^n=-1$.
		Faithfulness of $t_n$ implies that every $g\in G_{K_i}$ fixes
		$\xi_{2n}$; hence $\mu_{2n}\subset K_i$.
		Likewise, $G_L$ fixes $h$ and $\phi(h)=\eta$, so evaluation
		using $[(h,\sqrt[n]{\eta})]$ gives
		\[
		0=c(g)(h)
		=t_n\!\left(\frac{g^{-1}(\sqrt[n]{\eta})}{\sqrt[n]{\eta}}\right)
		\qquad(g\in G_L).
		\]
		Hence $\sqrt[n]{\eta}\in L$, and so $n\in\mathcal I_{a,b,c,m}$.
		
		Conversely, assume $n\in\mathcal I_{a,b,c,m}$.
		Our chosen roots satisfy $\xi_{2n}\in K_1\cap K_2\cap K_3$
		and $\sqrt[n]{\eta}\in L$.
		The lifts of the basis $e_1,e_2,e_3,h$, together with
		$\mu_n\subset k$, therefore show that $G_L$ fixes $T_n$.
		Hence $c|_{G_L}=0$, so $\partial(1/n)$ restricts to zero over $L$.
		Since $G_L$ also fixes $\mathcal A$, inflation--restriction gives
		\[
		0\longrightarrow H^1(G,\mathcal A)
		\xrightarrow{\mathrm{inf}}H^1(k,\mathcal A)
		\xrightarrow{\mathrm{res}}H^1(L,\mathcal A),
		\qquad G=\Gal(L/k).
		\]
		It follows that $c$ descends to a cocycle on $G$, and it suffices to prove
		$[c]=0$ in $H^1(G,\mathcal A)$.
		
		Put $K_4=k(\sqrt{m\Delta})$, $H_i=\Gal(L/K_i)$ and
		$C_i=\QQ/\ZZ(\chi_i)$ for $1\le i\le4$.
		By Proposition~\ref{prop:alg-d4} and
		Lemma~\ref{lem:alg-dual-basis},
		$\Lambda_0:=\bigoplus_i\ZZ e_i\subset\Lambda$ has quotient
		$\ZZ/2$, generated by $h$, with $2h=\sum_i e_i$.
		Dualizing with $\QQ/\ZZ$ gives the exact sequence of $G$-modules
		\[
		0\longrightarrow\ZZ/2\xrightarrow{\iota}\mathcal A
		\xrightarrow{\pi}\bigoplus_i C_i\longrightarrow0,
		\qquad \pi(\ell)=(\ell(e_i))_i.
		\]
		Here $\iota(1)(e_j)=0$ and $\iota(1)(h)=1/2$.
		Choose $\ell_i\in\mathcal A$ with
		$\ell_i(e_j)=\delta_{ij}/2$ and $\ell_i(h)=1/4$.
		These values respect $2h=\sum_j e_j$.
		The same relation and $g(e_j)=(-1)^{\chi_j(g)}e_j$ give
		$g^{-1}h=h-\sum_j\chi_j(g)e_j$, hence
		\[
		(g\ell_i-\ell_i)(e_j)=0,\qquad
		(g\ell_i-\ell_i)(h)=-\frac{\chi_i(g)}2
		\quad\text{in }\QQ/\ZZ.
		\]
		We obtain $g\ell_i-\ell_i=\iota(\chi_i(g))$ and
		$\iota_*(\chi_i)=0\in H^1(G,\mathcal A)$.
		Since $L$ is the compositum of the quadratic extensions defined by
		the $\chi_i$, these characters span
		$H^1(G,\ZZ/2)=\Hom(G,\ZZ/2)$.
		The cohomology sequence therefore gives an injection
		\begin{equation}\label{eq:trans-H1-detection}
			H^1(G,\mathcal A)\lhook\joinrel\longrightarrow
			\prod_i H^1(G,C_i).
		\end{equation}
		
		For each $i$, restriction
		$H^1(G,C_i)\to H^1(H_i,\QQ/\ZZ)$ is injective:
		if $\chi_i=0$ it is the identity, and otherwise
		inflation--restriction and cyclic periodicity give
		\[
		H^1(G/H_i,C_i)
		=\widehat H^{-1}(G/H_i,C_i)
		=\frac{\QQ/\ZZ}{2(\QQ/\ZZ)}=0.
		\]
		For $i=1,2,3$, the lift $[(e_i,\xi_{2n})]$ is $H_i$-fixed,
		so the computation \eqref{eq:trans-coordinate-cocycle} gives
		$\pi_{i,*}[c]=0$.
		
		For $i=4$, it suffices to find $s_4\in K_4^*$ with
		$s_4^n=\phi(e_4)=-\eta^2$: then $[(e_4,s_4)]$ is $H_4$-fixed
		and the same evaluation gives $\pi_{4,*}[c]=0$.
		Since $rs=q_m/(abc)$, direct expansion gives
		\[
		\eta^2=\frac{r^2+s^2-2rs}{4}
		=\frac{m+\Delta-2q_m}{4abc}\in k(q_m)=K_4.
		\]
		
		If $\xi_{2n}\in K_4$, then for every $g\in H_4$,
		\[
		\alpha_g:=\frac{g(\sqrt[n]{\eta})}{\sqrt[n]{\eta}}
		\in\mu_{2n}\subset K_4,
		\qquad \alpha_g^n=\frac{g(\eta)}{\eta}\in\{1,-1\}.
		\]
		Since $g(\alpha_g)=\alpha_g$ and $g^2=1$ in $G$, we get
		\[
		\alpha_g^2=g(\alpha_g)\alpha_g
		=\frac{g^2(\sqrt[n]{\eta})}{\sqrt[n]{\eta}}=1.
		\]
		Thus $\bigl(\sqrt[n]{\eta}\bigr)^2$ is fixed by $H_4$ and
		belongs to $L^{H_4}=K_4$.  We may take
		\[
		s_4:=\xi_{2n}\bigl(\sqrt[n]{\eta}\bigr)^2\in K_4^*,
		\qquad s_4^n=\xi_{2n}^n\eta^2=-\eta^2.
		\]
		
		If $\xi_{2n}\notin K_4$, then $\xi_{2n}\notin k$;
		since $\xi_{2n}\in K_i$ and $[K_i:k]\le2$,
		$K_1=K_2=K_3=k(\xi_{2n})$.
		It follows that $L=K_4(\xi_{2n})$ is quadratic over $K_4$.
		Its nontrivial involution $\tau$ negates $q_a,q_b,q_c$ and
		fixes $q_m$.  The definitions of $r,s$ therefore give
		$\tau(r)=-r$ and $\tau(s)=-s$, hence $\tau(\eta)=-\eta$.
		We may therefore take
		\[
		s_4=N_{L/K_4}\!\left(\sqrt[n]{\eta}\right),\qquad
		s_4^n=N_{L/K_4}(\eta)=-\eta^2.
		\]
		All four coordinate classes vanish, so
		\eqref{eq:trans-H1-detection} gives $[c]=0$ in
		$H^1(G,\mathcal A)$, and hence $\partial(1/n)=0$.
		
		We have shown that $I$ contains an element of order $n$ exactly when
		$n\in\mathcal I_{a,b,c,m}$.
		The uniqueness of finite cyclic subgroups of $\QQ/\ZZ$ yields
		\eqref{trans-eq:quotient} and implies closure under divisors and least
		common multiples.
		Changing the square roots replaces $\eta$ by $\pm\eta^{\pm1}$;
		since $\mu_{2n}\subset L$ implies $-1\in L^{\times n}$,
		the criterion is independent of these choices.
	\end{proof}
	
	\begin{cor}\label{trans-cor:real}
		Let $k=\mathbf R$, let $a,b,c\in\mathbf R^*$, and suppose
		$m\ne0,4abc$.  With $\Delta=m-4abc$, one has
		\begin{equation}\label{eq:trans-real}
			\frac{\Br(U_m)}{\Br_1(U_m)}
			\simeq
			\begin{cases}
				\ZZ/2,&a\Delta<0,\ b\Delta<0,\ c\Delta<0,\\
				0,&\text{otherwise}.
			\end{cases}
		\end{equation}
	\end{cor}
	
	\begin{proof}
		Since $\mu_n\subset\mathbf R$ only for $n=1,2$,
		Theorem~\ref{trans-thm:explicit} shows that only $n=2$ can contribute
		a nontrivial quotient.
		For this order, $\mu_4\subset K_i$ is equivalent to
		$K_i=\mathbf C$, hence to $a\Delta<0$, $b\Delta<0$,
		$c\Delta<0$, respectively.  If all three inequalities hold,
		then $L=\mathbf C$ and $\eta$ is automatically a square in $L$.
		These sign conditions are therefore both necessary and sufficient.
	\end{proof}
	
	\begin{cor}\label{trans-cor:odd-number-field}
		If $k$ is a number field, we write $w_F=|\mu(F)|$ for each number
		field $F$, and put
		\[
		E=\gcd\left(w_k,\frac{w_{K_1}}2,
		\frac{w_{K_2}}2,\frac{w_{K_3}}2\right),
		\qquad h=[L(\sqrt[E]{\eta}):L].
		\]
		Then $h\mid E$ and
		\begin{equation}\label{eq:trans-number-field-degree}
			\boxed{\displaystyle
				\frac{\Br(U_m)}{\Br_1(U_m)}\simeq\ZZ/N,
				\qquad N=\frac{E}{h}.}
		\end{equation}
		In particular, $N=\max\mathcal I_{a,b,c,m}$ and $N\mid w_k$.
		If $E=1$,  the transcendental part is trivial. For $E=2$ case, it is
		$\ZZ/2$ exactly when $\eta$ is a square in $L$.
	\end{cor}
	
	\begin{proof}		
		Over a number field, the cyclotomic conditions in
		Theorem~\ref{trans-thm:explicit} are exactly $n\mid E$:
		they say $n\mid w_k$ and $2n\mid w_{K_i}$ for every $i$.
		Choose $u$ with $u^E=\eta$.  Since $\mu_E\subset k\subset L$,
		the extension $L(u)/L$ is Galois, and
		\[
		\Gal(L(u)/L)\hookrightarrow\mu_E,\qquad
		g\longmapsto g(u)/u
		\]
		identifies $\Gal(L(u)/L)$ with $\mu_h$, where $h=[L(u):L]\mid E$.
		For every divisor $n$ of $E$,
		\[
		\eta\in L^{\times n}
		\quad\Longleftrightarrow\quad u^{E/n}\in L
		\quad\Longleftrightarrow\quad h\mid E/n.
		\]
		The admissible indices in Theorem~\ref{trans-thm:explicit}
		are therefore precisely the divisors of $E/h$.
		Their largest element gives \eqref{eq:trans-number-field-degree}.
	\end{proof}
	
	\begin{rmk}\label{rmk:trans-odd-part}
		For odd $n$, the condition $\mu_n\subset k$ implies
		$\mu_{2n}\subset k$. Hence
		\[
		\Br(U_m)/\Br_1(U_m)
		\text{ contains a subgroup of order }n
		\]
		if and only if
		\[
		\mu_n\subset k
		\qquad\text{and}\qquad
		\eta\in L^{\times n}.
		\]
		For odd $n$, the latter condition can be tested in
		$K_4=k(\sqrt{m\Delta})$:
		\[
		\eta\in L^{\times n}
		\quad\Longleftrightarrow\quad
		\eta^2\in K_4^{\times n}.
		\]
		Consequently, for a number field $k$, the subgroup of elements of odd
		order in the transcendental part is cyclic of order
		\[
		N_{\mathrm{odd}}
		=
		\max\{n\mid w_k:\ n\text{ odd and }
		\eta^2\in K_4^{\times n}\}.
		\]
		In particular, it is trivial when $\mu(k)=\{\pm1\}$.
		
	\end{rmk}
	
	\section{Consequences and examples}\label{sec:consequences}
	
	\subsection{Fields with a real embedding}
	
	Let $k$ be a number field with a real embedding
	$v_0:k\hookrightarrow\mathbf R$ (for example, $k=\QQ$).
	Every root of unity in $k$ maps to a real root of unity, so
	$\mu(k)=\{\pm1\}$ and $-1\notin k^{\times2}$.
	Theorem~\ref{trans-thm:explicit} therefore shows that
	$\Br(U_m)/\Br_1(U_m)$ is either $0$ or $\ZZ/2$.
	If the transcendental part is nonzero, the cyclotomic conditions force
	\[
	K_1=K_2=K_3=k(\sqrt{-1}),\qquad
	-a\Delta,-b\Delta,-c\Delta\in k^{\times2}.
	\]
	These square conditions imply that $ab,ac,bc$ are squares in $k$
	and that
	\[
	-D=\frac{-a\Delta}{(ab)(ac)}\in k^{\times2}.
	\]
	Choose $d\in k$ with $d^2=-D$ and $v_0(d)>0$.
	The normalization \eqref{eq:alg-coordinate-change} is then defined
	over $k$ and identifies $U_m$ with the ordinary Markoff surface
	\[
	U_m\simeq_k\{X^2+Y^2+Z^2-XYZ=4-d^2\}.
	\]
 	All square tests and the Galois-module calculation below take place
    over $k$. The real embedding ensures that $v_0(p)>0$ for $p=1+d/2$,
    so $-p\notin k^{\times2}$. This gives the following description of
    the full Brauer group over such number fields.
	
	\begin{thm}\label{thm:full-real-number-fields}
		Let $k$ be a number field with a real embedding
		$v_0:k\hookrightarrow\mathbf R$, let $a,b,c\in k^*$, and let
		$m\in k\setminus\{0,4abc\}$.
		Put $\Delta=m-4abc$, $M=m/(abc)$, $D=M-4$, and let
		$\chi_1,\ldots,\chi_4$ be the characters of
		$a\Delta,b\Delta,c\Delta,m\Delta$.  Set
		\[
		N=\#\{i:\chi_i\ne0\},\qquad
		\rho=\dim_{\FF_2}\langle\chi_1,\ldots,\chi_4\rangle.
		\]
		When $-a\Delta,-b\Delta,-c\Delta$ are all squares in $k$,
		we have $-D\in k^{\times2}$; in this case choose $d\in k$
		with $d^2=-D$ and $v_0(d)>0$, and put
		\[
		p=1+\frac d2,\qquad q=1-\frac d2,\qquad M=4pq.
		\]
		Then the full Brauer group is as follows.
		\begin{enumerate}
			\item If at least one of $-a\Delta,-b\Delta,-c\Delta$ is not a
			square in $k$, or if all three are squares but none of
			$p,q,-q$ is a square in $k$, then
			\[
			\Br(U_m)=\Br_1(U_m)\simeq\Br(k)\oplus
			\begin{cases}
				(\ZZ/2)^{N-\rho},&N<4,\\
				(\ZZ/2)^{5-\rho},&N=4,\ \text{no odd-weight relation},\\
				\ZZ/4\oplus(\ZZ/2)^{3-\rho},
				&N=4,\ \text{an odd-weight relation}.
			\end{cases}
			\]
			Here an odd-weight relation is a relation
			$\sum_i\epsilon_i\chi_i=0$ with
			$\sum_i\epsilon_i=1$ in $\FF_2$.
			\item Suppose $-a\Delta,-b\Delta,-c\Delta$ are all squares in $k$
			and at least one of $p,q,-q$ is a square in $k$.  Then
			\[
			\Br(U_m)/\Br_1(U_m)\simeq\ZZ/2,\qquad
			\Br_1(U_m)\simeq\Br(k)\oplus(\ZZ/2)^\ell,
			\]
			where
			\[
			\ell=
			\begin{cases}
				4,&M\in k^{\times2},\\
				2,&-M\in k^{\times2},\\
				3,&\text{otherwise}.
			\end{cases}
			\]
			The full group is given by the following two cases:
			\begin{subequations}\label{eq:full-real-number-fields}
				\begin{enumerate}[label={\upshape(\alph*)}]
					\item If $-q\in k^{\times2}$ and $p\notin k^{\times2}$,
					then $\ell=3$ and the extension is nonsplit, with
					\begin{equation}
						\Br(U_m)/\Br(k)\simeq\ZZ/4\oplus(\ZZ/2)^2.
					\end{equation}
					\item If $p\in k^{\times2}$ or $q\in k^{\times2}$,
					the extension splits, with
					\begin{equation}
						\Br(U_m)/\Br(k)\simeq(\ZZ/2)^{\ell+1}.
					\end{equation}
				\end{enumerate}
			\end{subequations}
		\end{enumerate}
		In all cases, the constant subgroup splits off:
		\[
		\Br(U_m)\simeq\Br(k)\oplus\bigl(\Br(U_m)/\Br(k)\bigr).
		\]
	\end{thm}
	
	\begin{proof}
		Put $A=\Br_1(U_m)/\Br(k)$ and $B=\Br(U_m)/\Br(k)$.
		The formula for the algebraic part in (1) follows from
		Theorem~\ref{thm:alg-brauer}. If one of $-a\Delta,-b\Delta,-c\Delta$
		is not a square, Theorem~\ref{trans-thm:explicit} gives $B=A$.
		 Assume henceforth that
		$-a\Delta,-b\Delta,-c\Delta$ are squares.
		The normalization \eqref{eq:alg-coordinate-change} gives the ordinary Markoff
		surface with parameter $M=4-d^2$ and characters
		$[-1],[-1],[-1],[-M]$.  Theorem~\ref{thm:alg-brauer}
		now gives $A=(\ZZ/2)^\ell$ as stated.  Put
		\[
		t=\sqrt M,\qquad F=k(t),\qquad L=F(i),\qquad
		\eta=\frac{t-id}{2},\qquad
		\kappa=\frac{(1-i)(d+2)+(1+i)t}{4}.
		\]
		Using $t^2-(d+2)^2=-2d(d+2)$, we complete the square:
		\[
		\begin{aligned}
			p\eta
			&=\frac{(d+2)(t-id)}4 
			=\frac{i}{8}\bigl(t^2-(d+2)^2-2i(d+2)t\bigr)\\
			&=\frac{i}{8}\bigl(t-i(d+2)\bigr)^2 
			=\left(\frac{(1+i)\bigl(t-i(d+2)\bigr)}4\right)^2\\
			&=\left(\frac{(1-i)(d+2)+(1+i)t}{4}\right)^2=\kappa^2.
		\end{aligned}
		\]
		Since $\kappa^2=p\eta$, Theorem~\ref{trans-thm:explicit} gives
		\[
		B/A\simeq\ZZ/2
		\iff\eta\in L^{\times2}
		\iff p\in L^{\times2}
		\iff\{p,q,-q\}\cap k^{\times2}\ne\varnothing.
		\]
		For the last equivalence, the quadratic square test gives
		\[
		\begin{aligned}
			p\in L^{\times2}
			&\iff p\in k(i)^{\times2}\text{ or }p/M\in k(i)^{\times2}\\
			&\iff p\in k(i)^{\times2}\text{ or }q\in k(i)^{\times2},
		\end{aligned}
		\]
		where the second line uses $p/M=1/(4q)$.
		The first line follows by squaring $x+y\sqrt M$ and comparing
		coefficients: if $\sqrt M\notin k(i)$, the cross term forces
		$xy=0$; otherwise the assertion is immediate.
		Likewise, a nonzero element of $k$ is a square in $k(i)$
		exactly when it or its negative is a square in $k$, since
		$(x+iy)^2\in k$ with $x,y\in k$ forces $xy=0$.
		As $v_0(p)>0$, $-p\notin k^{\times2}$, leaving precisely
		$p,q,-q$.
		
		Suppose $u\in L$ and $u^2=\eta$.
		The extension $0\to A\to B\to\ZZ/2\to0$ splits exactly when
		the nonzero quotient class has a lift of order $2$.
		By Theorem~\ref{thm:complete-affine-brauer} and
		\eqref{eq:brauer-restriction-duality}, this means finding
		$f\in\Hom_G(T_2,\FF_2)$ with $f(\tau)=1$, where
		$G=\Gal(L/k)$ and
		\[
		\tau=[(0,-1)],\qquad v_j=[(e_j,i)]\ (1\le j\le3),\qquad
		v_h=[(h,u)].
		\]
		These form a basis of $T_2$.
		Let $\psi$ and $\theta$ be the characters of $i$ and $t$, respectively,
		and put $E=e_1+e_2+e_3$.  The sign characters of
		$e_1,e_2,e_3,e_4$ are $\psi,\psi,\psi,\psi+\theta$.
		Since $e_4=2h-E$, taking $\theta(g)\in\{0,1\}$ gives
		\[
		\begin{aligned}
			g(h)&=\frac{(-1)^{\psi(g)}E+
				(-1)^{\psi(g)+\theta(g)}e_4}{2}\\
			&=(-1)^{\psi(g)+\theta(g)}h
			+(-1)^{\psi(g)}\theta(g)E
			\equiv h+\theta(g)E\pmod{2\Lambda}.
		\end{aligned}
		\]
		Moreover, $g(\tau)=\tau$ and $g(v_j)=v_j$:
		when $g(i)=-i$, one has
		$[(-e_j,-i)]=[(e_j,i)]$ because $(2e_j,-1)$ is a defining
		relation in $T_2$; otherwise $v_j$ is fixed directly.
		Under the projection $T_2\to\Lambda/2\Lambda$ in
		\eqref{eq:Tn-extension}, the difference
		$g(v_h)-v_h-\theta(g)(v_1+v_2+v_3)$ therefore maps to zero.
		The kernel of this projection is $\mu_2=\langle\tau\rangle$, so
		\begin{equation}\label{eq:real-T2-action}
			g(v_h)=v_h+\theta(g)(v_1+v_2+v_3)+c(g)\tau,
		\end{equation}
		for a unique $c(g)\in\FF_2$.
		Since $\tau,v_1,v_2,v_3$ are fixed, we have
		\[
		\begin{aligned}
			g_1g_2(v_h)
			&=g_1\bigl(v_h+\theta(g_2)(v_1+v_2+v_3)+c(g_2)\tau\bigr)\\
			&=g_1(v_h)+\theta(g_2)(v_1+v_2+v_3)+c(g_2)\tau\\
			&=v_h+\bigl(\theta(g_1)+\theta(g_2)\bigr)(v_1+v_2+v_3)
			+\bigl(c(g_1)+c(g_2)\bigr)\tau.
		\end{aligned}
		\]
		Comparing this with the formula for $g_1g_2(v_h)$, and using
		$\theta(g_1g_2)=\theta(g_1)+\theta(g_2)$, gives
		$c(g_1g_2)=c(g_1)+c(g_2)$ by uniqueness of the coefficient of
		$\tau$.  It follows that $c:G\to\FF_2$ is a character.
		Put $W=\langle\tau,v_1,v_2,v_3\rangle$, a trivial $G$-module.
		Dualizing $0\to W\to T_2\to\FF_2\to0$ and taking $G$-invariants gives
		\[
		\Hom_G(T_2,\FF_2)\longrightarrow\Hom(W,\FF_2)
		\xrightarrow{\delta}H^1(G,\FF_2),\qquad
		\delta(a,x_1,x_2,x_3)=a\,c+(x_1+x_2+x_3)\theta,
		\]
		where $a=f(\tau)$ and $x_j=f(v_j)$.
		The formula for $\delta$ follows by applying $f$ to $g(v_h)-v_h$ in
		\eqref{eq:real-T2-action}.
		A functional with $a=1$ therefore extends exactly when
		$c\in\langle\theta\rangle$.
		
		To evaluate this obstruction, put $z=\kappa/u$, so $z^2=p$.
		If $L\ne F=k(t)$, let $\sigma$ generate $\Gal(L/F)$.
		Then $\sigma(t)=t$, $\sigma(i)=-i$, $\theta(\sigma)=0$,
		and $\sigma(h)=-h$.  By the definitions of $v_h$ and $c$,
		\[
		\begin{aligned}
			c(\sigma)\tau
			&=\sigma(v_h)-v_h\\
			&=[(-2h,\sigma(u)/u)]\\
			&=[(0,\eta\,\sigma(u)/u)]
			=[(0,\sigma(u)u)].
		\end{aligned}
		\]
		Here the third line adds the defining relation
		$(2h,\phi(h))=(2h,\eta)$ from \eqref{eq:definition-Tn},
		and the last equality uses $u^2=\eta$.
		Since $\tau=[(0,-1)]$ and $\mu_2\hookrightarrow T_2$ is injective,
		this proves $(-1)^{c(\sigma)}=\sigma(u)u$.
		Direct expansion of the formula for $\kappa$ gives
		\[
		\sigma(\kappa)\kappa
		=\frac{(d+2+t)^2+(t-d-2)^2}{16}
		=\frac{(d+2)^2+t^2}{8}
		=\frac{d+2}{2}=p,
		\]
		where $t^2+d^2=4$.  Finally, $u=\kappa/z$ and $z^2=p$ give
		\[
		(-1)^{c(\sigma)}=\sigma(u)u
		=\frac{\sigma(\kappa)\kappa}{\sigma(z)z}
		=\frac{p}{\sigma(z)z}=\frac{z}{\sigma(z)}.
		\]
		As $\ker\theta=\Gal(L/F)$, this shows that
		$c\in\langle\theta\rangle$ exactly when $z\in F$,
		or equivalently $p\in F^{\times2}$.  If $L=F$, then $z\in F$,
		and $\ker\theta=1$ implies $c\in\langle\theta\rangle$.
		By the square test in $F=k(\sqrt M)$,
		the extension therefore splits exactly when $p$ or $q$ is a
		square in $k$, and then $B=A\oplus\ZZ/2$, giving (2b).
		If it does not split, the nontrivial quotient forces $-q$ to be
		a square and $p$ not to be a square.  It follows that $v_0(q)<0$
		and $v_0(M)<0$, while $-M=4p(-q)$ is not a square,
		so $\ell=3$.
		Since $2A=0$, the nonsplit extension is
		$B\simeq\ZZ/4\oplus(\ZZ/2)^{\ell-1}=\ZZ/4\oplus(\ZZ/2)^2$,
		giving (2a).  If the quotient is trivial, $B=A$.
		Finally, the constants split by
		Theorem~\ref{thm:complete-affine-brauer}.
	\end{proof}

	\subsection{The ordinary Markoff specialization}\label{subsec:ordinary}
	
	Suppose $ab$ and $ac$ are squares in $k$.  Choose $r_1^2=bc$,
	$r_2^2=ac$ and $r_3=r_1r_2/c$, so that $r_3^2=ab$ and
	$r_1r_2r_3=abc$.  The substitution
	\[
	X=x/r_1,\qquad Y=y/r_2,\qquad Z=z/r_3
	\]
	is defined over $k$ and identifies $U_m$ with the ordinary Markoff
	surface of parameter $M=m/(abc)$.  Put $D=M-4$ and
	$K=k(\sqrt M,\sqrt D)$, and choose
	$\eta=(\sqrt M-\sqrt D)/2$ as before.  The algebraic part is
	\begin{equation}\label{eq:standard-algebraic}
		\Br_1(U_m)/\Br(k)\simeq
		\begin{cases}
			0,&D\in k^{\times2},\\
			(\ZZ/2)^2,&D\notin k^{\times2},\ MD\in k^{\times2},\\
			(\ZZ/2)^4,&D\notin k^{\times2},\ M\in k^{\times2},\\
			(\ZZ/2)^3,&[M],[D]\text{ are independent}.
		\end{cases}
	\end{equation}
	This recovers the formulas of \cite[Theorem~3.4]{CTWX20}.
	
	For the transcendental part, define
	\begin{equation}\label{eq:ordinary-I}
		\mathcal I_{\mathrm{ord}}=\{n\ge1:\mu_n\subset k\text{ and }
		-1,\eta\in K^{\times n}\}.
	\end{equation}
	Theorem~4.8 of \cite{CTWX20} gives
	\begin{equation}\label{eq:ordinary-transcendental}
		\Br(U_m)/\Br_1(U_m)
		\simeq\varinjlim_{n\in\mathcal I_{\mathrm{ord}}}\ZZ/n.
	\end{equation}
	The transition maps are the natural embeddings in $\QQ/\ZZ$.
	The sign difference between the root used there and our $\eta$
	does not change \eqref{eq:ordinary-I}, because $-1$ is also required
	to be an $n$-th power.  Over a number field, the group in
	\eqref{eq:ordinary-transcendental} is finite cyclic.
	Formula~\eqref{eq:full-brauer-finite} also determines the extension
	of this cyclic group by \eqref{eq:standard-algebraic}.
	
	We verify explicitly that the conditions in \eqref{eq:ordinary-I} agree with
	those in Theorem~\ref{trans-thm:explicit}.  Put $F=k(\sqrt D)$.
	The square-class assumptions give $K_1=K_2=K_3=F$ and $L=K$.
	If $\mu_{2n}\subset F$, then $-1\in K^{\times n}$.
	Conversely, assume $\mu_n\subset k$ and $-1,\eta\in K^{\times n}$.
	If $K=F$, the power condition on $-1$ supplies a primitive
	$2n$-th root of unity in $F$ (for odd $n$ this already follows
	from $\mu_n\subset k$).  If $[K:F]=2$, the nontrivial automorphism
	changes the sign of $\sqrt M$ and fixes $\sqrt D$, so
	$N_{K/F}(\eta)=-1$.  Taking the norm of an $n$-th root of $\eta$
	shows that $-1\in F^{\times n}$ and gives the same cyclotomic
	containment.  Hence the two index sets coincide.
	
	\begin{example}\label{ex:mixed-brauer}
		Let $k=\QQ$.  Theorem~\ref{thm:full-real-number-fields} gives
		\begin{enumerate}
			\item For $U_{-12}: x^2+y^2+z^2-xyz=-12$,
			\[
			\begin{aligned}
				\Br_1(U_{-12})&\simeq\Br(\QQ)\oplus(\ZZ/2)^3,\\
				\Br(U_{-12})&\simeq\Br(\QQ)\oplus\ZZ/4\oplus(\ZZ/2)^2.
			\end{aligned}
			\]
			\item For $U_7: x^2+y^2+4z^2-xyz=7$,
			\[
			\begin{aligned}
				\Br_1(U_7)&\simeq\Br(\QQ)\oplus(\ZZ/2)^3,\\
				\Br(U_7)&\simeq\Br(\QQ)\oplus(\ZZ/2)^4.
			\end{aligned}
			\]
		\end{enumerate}
		
			In both cases, the algebraic part is $(\ZZ/2)^3$ and the
		transcendental part is $\ZZ/2$. However, the extension is nonsplit in the
		first case but split in the second, which means the algebraic and
		transcendental parts alone do not determine the full Brauer group.
	\end{example}

	\section{Integral Hasse failures over \texorpdfstring{$\QQ$}{Q}}
	\label{sec:integral-failures}
	
	In this section, the coefficients and the parameter are integers.
	Let $\mathcal U_m\subset\Aff^3_{\ZZ}$ denote the integral model
	$ax^2+by^2+cz^2-xyz=m$, with generic fibre $U_m$.  Write
	\[
	\mathcal U_m(\mathbb A_{\ZZ})
	=U_m(\mathbf R)\times\prod_p\mathcal U_m(\ZZ_p).
	\]
	Its integral Brauer--Manin set consists of tuples $(P_v)_v$ satisfying
	\[
	\sum_v\operatorname{inv}_v B(P_v)=0
	\qquad \forall B\in\Br(U_m).
	\]
	The integral Hasse failure means that this product of local sets is
	nonempty while $\mathcal U_m(\ZZ)$ is empty.  In the following two
	families, the failure is explained by a single quaternion class in $\Br(X_m)$.
	Local solubility follows by lifting smooth points modulo each prime
	to $\ZZ_p$, giving solutions to every congruence as well.
	
	\subsection{A positive family}
	
	The norm identities used for ordinary Markoff surfaces in
	\cite{GS22,CTWX20} and for their one-coefficient twists in
	\cite{SChen21} also give explicit families with three distinct
	coefficient square classes.  We give the local calculations directly.
	\begin{prop}\label{integral-prop:positive-family}
		Let $a,b,c$ be positive pairwise coprime odd integers satisfying
		\[
		a\equiv b\equiv c\equiv1\pmod8,
		\qquad ac\equiv-1\pmod3.
		\]
		Assume that every prime divisor of $abc$ is congruent to
		$\pm1$ modulo $8$.  Let $d$ be a positive odd integer coprime
		to $abc$, all of whose prime divisors are congruent to
		$\pm1$ modulo $8$, and let $m=4abc+2bd^2$.  Then
		\[
		\mathcal U_m(\mathbb A_{\ZZ})\ne\varnothing,
		\qquad
		\mathcal U_m(\mathbb A_{\ZZ})^{\Br(U_m)}=\varnothing.
		\]
		In particular, the equation has no solution in $\ZZ^3$.
		The obstruction is given by $B=(x^2-4bc,2)$.
	\end{prop}
	
	\begin{proof}
		Put $f=ax^2+by^2+cz^2-xyz-m$.
		The identity
		\begin{equation}\label{integral-eq:positive-norm}
			(2by-xz)^2-8b^2d^2
			=(x^2-4bc)(z^2-4ab)
		\end{equation}
		will be used to compute the local invariants.  Since
		$b\Delta=2b^2d^2$, Theorem~\ref{proj:classification} gives
		$B\in\Br(X_m)\subset\Br_1(U_m)$.
		
		We first check local solubility.  At $2$, the point $(1,1,1)$
		solves the equation modulo $2$, and the derivative with respect
		to $y$ is odd.  At $3$, we have $3\nmid abcd$,
		and $m\equiv b(ac-1)\not\equiv0\pmod3$.  The nonsingular
		binary form $by^2+cz^2$ represents this nonzero residue, giving
		a smooth point with $x=0$.  A nonsingular binary quadratic form
		over an odd finite field represents every nonzero element:
		this follows either from its split form or from surjectivity
		of the norm from the quadratic field extension.
		
		Let $p\ge5$.  If $p\mid a$, take $x=0$ and represent
		$m\equiv2bd^2\ne0$ by $by^2+cz^2$; if $p\mid c$, take
		$z=0$ and use $ax^2+by^2$.  Pairwise coprimality ensures these
		forms are nonsingular.  If $p\mid b$, the point
		$(1,a+c,1)$ solves the equation modulo $p$ and has
		$\partial f/\partial y=-1$.  For $p\nmid abc$, use $x=0$
		when $p\nmid m$.  When $p\mid m$, choose $t\in\FF_p^*$
		with $t^2\ne4bc$; such a choice exists since at most two
		nonzero elements are excluded.  The form
		$by^2+cz^2-tyz$ is nonsingular and represents
		$-at^2\ne0$, giving a smooth point with $x=t$.
		Hensel's lemma lifts all these points.  A real point is
		$(0,\sqrt{m/b},0)$.
		
		We now compute the invariants.  At every odd prime dividing
		$abc d$, $2$ is a square and $B$ evaluates to zero.  At any other
		odd prime for which $2$ is not a square, \eqref{integral-eq:positive-norm}
		shows that $x^2-4bc$ cannot vanish modulo $p$, since $b,d$
		are units.  Its value is therefore a unit and a norm from the unramified
		quadratic extension, and again the invariant is zero.
		The real invariant is zero as well.
		
		At $2$, we have $m\equiv6\pmod8$.  If $x$ is odd, then
		$x^2-4bc\equiv5\pmod8$.  If $x$ is even, the equation modulo
		$2$ forces $y$ and $z$ to be both even or odd. The former is impossible
		modulo $4$, and the equation modulo $8$ then forces
		$x\equiv4\pmod8$.  It follows that
		\[
		(x^2-4bc)/4\equiv3\pmod8.
		\]
		Both cases give $\operatorname{inv}_2B=1/2$, by the formula
		$(u,2)_2=(-1)^{(u^2-1)/8}$ for odd units $u$.
		Hence the sum of invariants is always $1/2$.  The integral
		Brauer--Manin set is therefore empty; reciprocity then rules out
		integral solutions.
	\end{proof}
	
	\subsection{A negative family of coefficient twists}
	\label{subsec:integral-negative}
	
	The negative family of Ghosh and Sarnak
	\cite[Proposition~8.1(i)]{GS22} admits the following coefficient twist.
	Its obstruction is detected by a class on the projective surface.
	
	\begin{prop}\label{thm:integral-negative}
		Let $A$ be a positive odd integer and let $v$ be a positive odd
		integer such that $v^2>2A$ and every prime divisor of $v$ is
		congruent to $1$ or $3$ modulo $8$.  Put $m=4A-2v^2$ and consider
		the integral model
		\[
		\mathcal U_{A,v}:\quad Ax^2+y^2+z^2-xyz=m
		\subset\Aff^3_{\ZZ},
		\]
		with generic fibre $U_{A,v}$ and smooth projective completion
		$X_{A,v}$.  Then $\mathcal U_{A,v}(\ZZ_p)\ne\varnothing$ for every
		prime $p$ and $U_{A,v}(\mathbb R)\ne\varnothing$, but
		\[
		\mathcal U_{A,v}(\mathbb A_{\ZZ})^{\Br(U_{A,v})}
		=\varnothing,\qquad
		\mathcal U_{A,v}(\ZZ)=\varnothing.
		\]
		More precisely, the class
		\[
		B=(x^2-4,-2)\in\Br(X_{A,v})
		\]
		has invariant $1/2$ at every integral $2$-adic point and invariant
		zero at all other integral local points, including the real place.
	\end{prop}
	
	\begin{proof}
		Here $\Delta=m-4A=-2v^2$, so $B$ is the quaternion class in
        $\Br(X_{A,v})$ given by Theorem~\ref{proj:classification}.  
        The conditions $m<0$ and $\Delta\ne0$ ensure that the surfaces
        are smooth.
		
		We first establish local solubility.  At $2$, put $y=z=1$ and solve
		\[
		Ax^2-x+2=m.
		\]
		The value at $x=0$ is even, while the derivative $2Ax-1$ is odd.
		Hensel's lemma therefore gives a solution in $\ZZ_2$.
		At $3$, we have $m\equiv A+v^2\pmod3$.  If $3\nmid m$, set $x=0$
		and represent the nonzero residue $m$ by $y^2+z^2$ over $\FF_3$;
		the resulting point is smooth.  If $3\mid m$ and $A\equiv2\pmod3$,
		the point $(1,1,0)$ solves the equation modulo $3$ and the partial derivative
		about $y$ is $2$.  The only other possibility with $3\mid m$ is
		$3\mid A$ and $3\mid v$.  Then $(2,1,1)$ is a solution modulo $3$,
		and the partial derivative with respect to $x$, namely $4A-1$, is a unit.
		
		Let $p\ge5$.  If $p\nmid m$, a representation of $m$ by
		$y^2+z^2$ over $\FF_p$, with $x=0$, gives a smooth point.
		If $p\mid m$ and $p\nmid A$, put $x=1$ and represent $-A$ by
		\[
		y^2+z^2-yz
		\]
		over $\FF_p$.  This nondegenerate binary quadratic form represents
		every nonzero element of $\FF_p$, and the chosen point is smooth
		because its value is nonzero.  If $p\mid m$ and $p\mid A$, the
		point $(2,1,1)$ is a solution modulo $p$, with derivative
		$4A-1\equiv-1\pmod p$ in $x$.  All these points lift by Hensel's
		lemma.  Finally, the polynomial $(A+2)t^2-t^3$ takes the value
		$m<0$ for some positive real $t$, yielding a real point.
		
		For the invariant at $2$, note that $m\equiv2\pmod8$.
		If $x$ is odd, then $x^2-4\equiv5\pmod8$ and
		$(x^2-4,-2)_2=-1$.  Suppose that $x$ is even.  Parity forces $y,z$
		to have the same parity.  They cannot both be even: in that case
		the left side of the equation is divisible by $4$, whereas
		$m\equiv2\pmod4$.  Hence $y,z$ are odd.  If $x\equiv2\pmod4$,
		then $Ax^2+y^2+z^2\equiv6\pmod8$ and $xyz\equiv2$ or $6\pmod8$,
		which is incompatible with $m\equiv2\pmod8$.  Hence $4\mid x$.
		The equation then gives $xyz\equiv0\pmod8$, so $8\mid x$.
		It follows that
		\[
		\frac{x^2-4}{4}\equiv-1\pmod8,
		\qquad (x^2-4,-2)_2=(-1,-2)_2=-1.
		\]
		In both cases $\operatorname{inv}_2B=1/2$.
		
		At $3$, the element $-2$ is a square in $\QQ_3$, so $B$ vanishes.
		Let $p\ne2,3$.  If $-2$ is a square in $\QQ_p$, the same conclusion
		holds.  Otherwise, no integral local point can have
		$x\equiv\pm2\pmod p$: substituting this value into the equation
		would give
		\[
		(y\mp z)^2\equiv-2v^2\pmod p,
		\]
		contrary to the nonsquareness of $-2$.  Here $p\nmid v$, since
		all prime divisors of $v$ split in $\QQ(\sqrt{-2})$ by hypothesis.
		It follows that $x^2-4$ is a $p$-adic
		unit.  The Hilbert symbol of two units at an odd prime is trivial,
		so $\operatorname{inv}_pB=0$.
		
		At a real point we have $|x|>2$.  Indeed, if $|x|\le2$, then
		\[
		m=Ax^2+(y-xz/2)^2+(1-x^2/4)z^2\ge0,
		\]
		contrary to $m<0$.  Hence $x^2-4>0$ and the real invariant vanishes.
		The sum of the local invariants of $B$ is therefore always $1/2$.
		It follows that the integral Brauer--Manin set is empty and,
		by reciprocity, that no integral point exists.
	\end{proof}
	
	\begin{example}\label{ex:integral-negative-order-four}
		Taking $A=1$ in Proposition~\ref{thm:integral-negative} recovers
		\cite[Proposition~8.1(i)]{GS22}.  Taking $A=3$ and $v=3^n$ gives
		the integral Hasse failures
		\[
		3x^2+y^2+z^2-xyz=12-2\cdot9^n\qquad(n\ge1).
		\]
		For $n=1,2$, the character tuple is
		$([-6],[-2],[-2],[3])$, with $N=4$, $d=2$ and the odd-weight
		relation $\chi_1+\chi_2+\chi_4=0$.
		Theorem~\ref{thm:alg-brauer} gives
		\[
		\Br(U_{3,3^n})/\Br(\QQ)\simeq\ZZ/4\oplus\ZZ/2
		\qquad(n=1,2).
		\]
		This is the full group by Theorem~\ref{thm:full-real-number-fields},
		since the coefficient classes $[3],[1],[1]$ do not all agree.
	\end{example}
	
	\begin{rmk}\label{rmk:integral-hasse-beyond-brauer}	
		Integral Hasse failures need not be explained by the Brauer--Manin
		obstruction.  For the ordinary Markoff surface with
		$m=342=4+2\cdot13^2$, Ghosh--Sarnak
		\cite[Proposition~8.1(iii)]{GS22} prove that there are no integral
		points, whereas Colliot-Th\'el\`ene--Wei--Xu
		\cite[Proposition~5.11]{CTWX20} prove that the integral
		Brauer--Manin set is nonempty.  More generally, the same two
		conclusions hold for $m=4+2\ell^2$ when $\ell\ge13$ is prime,
		$\ell\equiv\pm4\pmod9$ and $\ell\equiv\pm3\pmod8$.
		There are infinitely many such primes, for example the primes in the
		progression $\ell\equiv13\pmod{72}$.  The proof of nonexistence
		uses reduction theory as well as reciprocity.
	\end{rmk}

\end{document}